\documentclass[12pt,reqno]{amsart}
\usepackage{amsmath,amssymb,epsfig,latexsym}

\DeclareMathSymbol{\Z}{\mathalpha}{AMSb}{"5A}

\newcommand{\df}{\dfrac}
\renewcommand{\i}{\infty}
\numberwithin{equation}{section}
 \theoremstyle{plain}
\newtheorem{theorem}{Theorem}[section]

\begin{document}
\title[Shifted and Shiftless Partition Identities II]
      {Shifted and Shiftless Partition Identities II}
\author{Frank G. Garvan}
\email{frank@math.ufl.edu}
\author{Hamza Yesilyurt}
\email{hamza@math.ufl.edu}
\address{Department of Mathematics\\
         University of Florida\\
         Gainesville, Florida 32611}
\begin{abstract}
Let $S$ and $T$ be sets of positive integers and let $a$ be a fixed positive
integer. An $a$-shifted partition identity has the form
$$
p(S,n)=p(T,n-a),\quad\mbox{for all $n \geq a$}.
$$
Here $p(S,n)$ is the number partitions of $n$ whose parts are
elements of $S$. For all known nontrivial shifted partition
identities, the sets $S$ and $T$ are unions of arithmetic
progressions modulo $M$ for some $M$. In 1987, Andrews found two
$1$-shifted examples ($M=32$, $40$) and asked whether there were any
more. In 1989, Kalvade responded with a further six. In 2000, the
first author found $59$   new $1$-shifted identities using a
computer search and showed how these could be proved using the
theory of modular functions.

Modular transformation of
certain shifted identities leads to shiftless partition identities.
Again let $a$ be a fixed positive integer, and $S$, $T$ be distinct
sets of positive integers. A shiftless partition identity has the form
$$
p(S,n)=p(T,n),\quad\text{for all $n \neq a$.}
$$

In this paper, we show, except in one case,  how all known
$1$-shifted and shiftless identities follow from a four parameter
theta function identity due to Jacobi. New shifted and shiftless
partition identities are proved.
\end{abstract}

\maketitle


\section{Introduction}

Let $S$ and $T$ be sets of positive integers. Let $a$ be a fixed
positive integer. An {\it $a$-shifted partition identity} has the
form
$$
p(S,n)=p(T,n-a),\qquad \mbox{for all $n\ge a$}.
$$

Andrews \cite{and} found the first nontrivial $1$-shifted
identities:
\begin{align*}
S& =\{n\, : \,\mbox{$n$ odd or }
    n\equiv \pm 4,\pm 6,\pm 8,\pm 10 \pmod{32}\},\\
T& =\{n\, : \, \mbox{$n$ odd or }
    n\equiv \pm 2,\pm 8,\pm 12,\pm 14\pmod{32}\};
\end{align*}
and
\begin{align*}
S&=\{n\, : \,n\equiv \pm 1,\pm 4,\pm 5,\pm 6,\pm 7,\pm 9,\pm 10,\pm 11,
                         \pm 13,\pm 15,
 \pm 16,\pm 19 \pmod{40}\}, \nonumber\\
T&=\{n\, : \, n\equiv \pm 1,\pm 3,\pm 4,\pm 5,\pm 9,\pm 10,\pm 11,\pm 14,
                        \pm 15,\pm 16,
 \pm 17,\pm 19 \pmod{40}\},\nonumber\\
\nonumber
\end{align*}
In these examples, each $S$ and $T$ is the union of arithmetic
progressions modulo $M$ for some $M$, namely $M=32$ and $M=40$. In
fact, each is a union of $24$ such arithmetic progressions. In 1989,
Kalvade \cite{KA} found six more identities with $M=42$, $48$ and
$60$, each also involving the union of $24$ arithmetic progressions.
In 1996, Alladi \cite{AG} considered more general shifted partition
identities in which parts in some residue classes are distinct. In
2000, the first author \cite{garvan00} found 59 new $1$-shifted
identities using a computer search and showed how such identities
could be proved using the theory of modular functions. One case
\cite[(4.1)]{garvan00} was proved in detail.


The generating function for $p(S,n)$ is an infinite product:
$$
\sum_{n\ge0} p(S,n) q^m = \prod_{n\in S} \frac{1}{(1-q^n)},
$$
where $|q|<1$.  We can write an $a$-shifted partition identity
as an equivalent $q$-series identity:
\begin{equation}
\prod_{n\in S} \frac{1}{(1-q^n)} - q^a \prod_{n\in T} \frac{1}{(1-q^n)} = 1.
\label{genid}
\end{equation}

The following two identities of Ramanujan which were later proved by
Rogers \cite{rogers} are well-known examples of shifted partition
identities
\begin{equation}
H(q)G(q^{11})-q^2G(q)H(q^{11})=1
\label{rr1}
\end{equation}
and
\begin{equation}
H(q^2)G(q^7)-qG(q^2)H(q^7)=
\prod_{n=1}^\infty \df{(1-q^{2n-1})}{(1-q^{14n-7})},
\label{rr2}
\end{equation}
where
\begin{equation*}
G(q) =
\prod_{n\equiv\pm1\pmod{5}} \frac{1}{(1-q^n)}
\;\;\text{and}\;\;
H(q) =
\prod_{n\equiv\pm2\pmod{5}} \frac{1}{(1-q^n)}
\;\;
\end{equation*}
are the Rogers-Ramanujan functions. Equation \eqref{rr1} is
equivalent to a $2$-shifted identity observed by Bressoud \cite{bre}.
See also Andrews \cite[(3)]{and}.

In \cite{garvan00}, the effect of modular transformations on shifted
identities was studied.
This led to certain shiftless partition identities.
As before, let $a$ be a fixed positive integer, and let $S$, $T$ be
distinct sets of positive integers. A {\it shiftless partition identity}
has the form:
$$
p(S,n) = p(T,n),\qquad  \mbox{for all $n\ne a$.}
$$
The simplest shiftless identity has $M=40$ and $a=2$. See
\cite[(1.5)]{garvan00} or Theorem \ref{40.1} below. Each shiftless
partition identity, that we find, can be written as an equivalent
$q$-series identity
\begin{equation}
\frac{1}{q^a} \prod_{n\in S} \frac{1}{(1-q^n)} - 
\frac{1}{q^a} \prod_{n\in T} \frac{1}{(1-q^n)} = 1.
\label{shiftid}
\end{equation}

In this paper we show, except for one case, that all known
$1$-shifted and shiftless partition identities follow from a certain
four parameter theta function identity (equation \eqref{jkb}) due to
Jacobi. The exceptional case which is given in Theorem \ref{72.2} is
similar to  the identities \eqref{rr1} and \eqref{rr2} above.

We note that in all identities considered the two products that occur on the 
left sides of
each of \eqref{genid} and \eqref{shiftid} can be written as a quotient
of an eta-product and a product of theta-functions. When we apply
a modular transformation 
$A=\begin{pmatrix} \alpha & \beta \\ \gamma  & \delta \end{pmatrix}
\in \Gamma_0(M)$, we obtain a shifted or shiftless identity.
Here $\gcd(\alpha,M)=1$, and the effect on $S$ and $T$ is to basically
multiply residue classes by $\alpha$ and reduce the products modulo $M$.
Thus we say two identities are equivalent under multiplication
by the group $U(M)=\left({\Z}/{M\Z}\right)^{\times}$
if one can be obtained from the other by applying a
modular transformation in $\Gamma_0(M)$. This is explained
in more detail with examples in \cite{garvan00}.
All shifted and shiftless identities considered can be grouped into
equivalence classes. 
One needs only prove one
identity in each class. Many of these identities are just special
cases of the four parameter Jacobi identity \eqref{jkb}. Others
follow from several theta function identities each of which is a
special case of \eqref{jkb}.

We find new shifted and shiftless partition identities which are not
equivalent to  $1$-shifted identities under modular transformation.
In any shifted or shiftless identity of given modulus $M$ we may
replace $q$ by $q^k$ and get another identity of modulus $kM$.
Hence, we may assume throughout that $\gcd(S \cup T)=1$.

\section{Theta function identities}
For $|q|<1$, let us define
\begin{equation*}
(a;q)_\i=(a)_\i:=\prod_{n=0}^{\i}(1-aq^n),
\end{equation*}
\begin{align*}
[a_1,\cdots a_k:q]&:=(a_1)_\i(q/a_1)_\i\cdots(a_k)_\i(q/a_k)_\i\;
,\text{and}\\ (a_1,\cdots
a_k:q)&:=(-a_1)_\i(-q/a_1)_\i\cdots(-a_k)_\i(-q/a_k)_\i.
\end{align*}



The main tool in our proofs is a four parameter $\theta-$function
identity due to Jacobi, namely
\begin{equation}\label{jkb}
[z,t,xty,zy/x:q]-[xt,zy,ty,z/x:q]=\df{z}{x}[y,x,xt/z,zty:q].
\end{equation}
The following equivalent form of \eqref{jkb} in terms of sigma
function is given in \cite[p.~451]{ww}
\begin{equation*}
\sum_{u,\,v,\,w}\sigma(z+u)\sigma(z-u)\sigma(v+w)\sigma(v-w)=0.
\end{equation*}
For convenience, we use the following reformulation of \eqref{jkb}
which in  most general form can also be found as an exercise in
\cite{ww}.
\begin{equation}\label{four}
-bc^2[b/c,a/x,a/y,xy/bc:q]+ac^2[a/c,b/x,b/y,xy/ac:q]=ab^2[a/b,c/x,c/y,xy/ab:q].
\end{equation}
In \eqref{four}, by replacing $a$ by $-a$ and $b$ by $-b$, we
obtain
\begin{equation}\label{four22}
-bc^2(b/c,a/x,a/y,xy/bc:q)+ac^2(a/c,b/x,b/y,xy/ac:q)=ab^2[a/b,c/x,c/y,xy/ab:q].
\end{equation}
We divide both sides of \eqref{four22} by the right hand side of
\eqref{four22} and use the trivial identity
\begin{equation}\label{tval}
(a:q):=\df{[a^2:q^2]}{[a:q]}=\df{[a^2:q^2]}{[a:q^2][aq:q^2]},
\end{equation}
to write each infinite product  with base $q^2$, we arrive at
\begin{align}
&\df{-c^2[b^2/c^2,a^2/x^2,a^2/y^2,x^2y^2/b^2c^2:q^2]}
{ab[b/c,bq/c,a/x,aq/x,a/y,aq/y,xy/bc,xyq/bc,a/b,aq/b,c/x,cq/x,c/y,cq/y,xy/ab,xyq/ab:q^2]}\notag\\
&+\df{c^2[a^2/c^2,b^2/x^2,b^2/y^2,x^2y^2/a^2c^2:q^2]}
{b^2[a/c,aq/c,b/x,bq/x,b/y,bq/y,xy/ac,xyq/ac,a/b,aq/b,c/x,cq/x,c/y,cq/y,xy/ab,xyq/ab:q^2]}\notag\\
&=1.\label{four2}
\end{align}

In applications of  \eqref{four2}, $a,\,b,\,c,\,x,\,y$ and $q$ will
be replaced by $q^a$,$q^b$,$q^c$,$q^x$,$q^y$ and $q^n$ for some
positive integer $n$ but for convenience we will refer to these
parameters simply as  $[a,b,c,x,y]$. Here we may assume that
$\gcd(a,b,c,x,y)=1$. Shifted partition identities are obtained when
each product in the numerator cancels with a product in the
denominator and the remaining products in the denominator of both
terms on the left side of \eqref{four2} are distinct. However, such
cancellation imposes a bound on $n$ and we found that \eqref{four2}
does not directly produce shifted partition identities for $n>41$.
Starting in section \ref{322}, for each modulus, we first list those
shifted  partition identities that are obtained directly from
\eqref{four2}. For most of the moduli, we have identities that do
not follow directly from \eqref{four2} but these identities can
still be proved through ``iteration" by using certain special cases
of \eqref{four2} together with \eqref{four}. In most cases, the
auxiliary identities that we use for iteration are themselves
shifted partition identities. For convenience, in our proofs, we
will use the shorthand notation $[a:n]$ for $[q^a:q^n]$ and $(a:n)$
for $(q^a:q^n)$.

We should remark that all shifted partition identities proved by
Andrews, Kalvade and Alladi are special cases of \eqref{four2}. In
fact, we were able to simplify the Macdonald identity that Kalvade
employed in her proofs and reduce it to the following special case
of \eqref{jkb}
\begin{equation}\label{kal}
(x,y:q)\big\{(x^2/y,xy^2:q)
-\df{x}{y}(y^2/x,x^2y:q)\big\}=[x^2,y^2,xy,x/y:q].
\end{equation}
The identity \eqref{kal} however has  many applications to shifted
partition identities other than those considered by Kalvade.
Replacing $q$ by $q^3$ and $y$ by $xq$, we obtain the famous
Quintuple Product Identity \cite[p.~80, Entry 28(iv)]{III}
\begin{equation}\label{qpp}
(x,xq,xq^2:q^3)\Bigr\{(x^3q:q^3)-x(x^{-3}q:q^3)\Bigl\}=[x^2,qx^2,q^2x^2,q:q^3].
\end{equation}
Some of our proofs use \eqref{qpp} alone in the form

\begin{align}
&[x^{-3}q,x^{-3}q^4,x^6q^2:q^6]-x[x^{3}q,x^{3}q^4,x^{-6}q^2:q^6]\notag\\
&=[q,q^2,x,xq,xq^2,xq^3,xq^4,xq^5,x^2q,x^2q^3,x^2q^5,x^3q,x^3q^4,x^{-3}q,x^{-3}q^4:q^6],\label{qp}
\end{align}
which we obtain  by expressing each infinite product in
\eqref{qpp} with base $q^6$ via \eqref{tval}.

\section{Modulus $M = 32$}\label{322}
There is only one  distinct class  of identities under
multiplication by the group $U(32)$. The identity in Theorem
\ref{32.1} below was first given by Andrews \cite{and}.

\begin{theorem}\label{32.1}

\begin{align*}
\hspace{0.5in}&\text{Let}\; S\equiv\pm\{1, 3, 4, 5, 6, 7, 8, 9, 10, 11, 13, 15\}\pmod {32},\; \text{and} \\
&T\equiv\pm\{1, 2, 3, 5, 7, 8, 9, 11, 12, 13, 14, 15\}\pmod{32}.\; \text{Then}\\
&p(S,n)=p(T,n-1),\; \text{for all}\; n\geq 1.
\end{align*}
\end{theorem}
\begin{proof}
This identity follows from \eqref{four2} with the choice of
parameters $[1, 2, 4, 12, 13]$; i.e.,
$a=q,\,b=q^2,\,c=q^4,\,x=q^{12},\,y=q^{13}$, and with $q$ replaced
by $q^{16}$.
\end{proof}

\section{Modulus $M = 40$}
There are two distinct classes  of identities under multiplication
by the group $U(40)$. Part (i) of Theorem \ref{40.1} below was also
given by Andrews \cite{and}.

\begin{theorem}\label{40.1}

\begin{align*}
\mathrm{(i)}\hspace{0.5in}&\text{Let}\; S\equiv\pm\{1, 4, 5, 6, 7, 9, 10, 11, 13, 15, 16, 19\}\pmod {40},\; \text{and} \\
&T\equiv\pm\{1, 3, 4, 5, 9, 10, 11, 14, 15, 16, 17, 19\} \pmod{40}.\; \text{Then}\\
&p(S,n)=p(T,n-1),\; \text{for all}\; n\geq 1.
\end{align*}
\begin{align*}
\mathrm{(ii)}\hspace{0.5in}&\text{Let}\; S\equiv\pm\{2, 3, 5, 7, 8, 9, 10, 11, 12, 13, 15, 17\}\pmod {40},\; \text{and} \\
&T\equiv\pm\{1, 3, 5, 7, 8, 10, 12, 13, 15, 17, 18, 19\}\pmod{40}.\; \text{Then}\\
&p(S,n)=p(T,n-2),\; \text{for all}\; n\geq 2.
\end{align*}
\end{theorem}
\begin{proof}
The identities (i) and (ii) follow from \eqref{four2} with the
choice of parameters $[1, 2, 5, 15, 16]$, $[1, 3, 4, 14, 16]$ and
with $q$ replaced by $q^{20}$ in both instances.
\end{proof}
\begin{theorem}\label{40.2}

\begin{align*}
\mathrm{(i)}\hspace{0.5in}&\text{Let}\; S\equiv\pm\{1, 3, 4, 5, 6, 7, 13, 15, 16, 17, 18, 19\}\pmod {40},\; \text{and} \\
&T\equiv\pm\{1, 2, 3, 5, 9, 11, 12, 15, 16, 17, 18, 19\} \pmod{40}.\; \text{Then}\\
&p(S,n)=p(T,n-1),\; \text{for all}\; n\geq 1.
\end{align*}
\begin{align*}
\mathrm{(ii)}\hspace{0.5in}&\text{Let}\; S\equiv\pm\{3, 4, 5, 6, 7, 8, 9, 11, 13, 14, 15, 17\}\pmod {40},\; \text{and} \\
&T\equiv\pm\{1, 3, 5, 8, 9, 11, 12, 14, 15, 17, 18, 19\} \pmod{40}.\; \text{Then}\\
&p(S,n)=p(T,n-3),\; \text{for all}\; n\geq 3.
\end{align*}
\begin{align*}
\mathrm{(iii)}\hspace{0.5in}&\text{Let}\; S\equiv\pm\{1, 3, 4, 5, 6, 7, 8, 13, 14, 15, 17, 19\}\pmod {40},\; \text{and} \\
&T\equiv\pm\{1, 2, 5, 6, 7, 8, 9, 11, 12, 13, 15, 19\} \pmod{40}.\; \text{Then}\\
&p(S,n)=p(T,n),\; \text{for all}\; n\neq 2.
\end{align*}
\begin{align*}
\mathrm{(iv)}\hspace{0.5in}&\text{Let}\; S\equiv\pm\{2, 3, 4, 5, 7, 9, 11, 13, 14, 15, 16, 17\}\pmod {40},\; \text{and} \\
&T\equiv\pm\{1, 2, 5, 7, 9, 11, 12, 13, 15, 16, 18, 19\} \pmod{40}.\; \text{Then}\\
&p(S,n)=p(T,n-2),\; \text{for all}\; n\geq 2.
\end{align*}
\end{theorem}
\begin{proof}
The identities (i)--(iv) follow from \eqref{four2} with the
following sets of parameters
\begin{equation*}
[1, 3, 2, 5, 7],\;[1, 5, 2, 7, 8],\;[1, 2, 4, 8, 9],\;[1, 3, 4,
12, 17]
\end{equation*}
and with $q$ replaced by $q^{20}$ in each instance.
\end{proof}

\section{Modulus $M = 42$}
There are three distinct classes  of identities under multiplication
by the group $U(42)$. The first identity in Theorem \ref{42.1} below
was first given  by Kalvade \cite{KA}.
Kalvade also conjectured the two identities, numbered (iv) and (vi),  in Theorem \ref{42.2}. \\\
\begin{theorem}\label{42.1}

\begin{align*}
\mathrm{(i)}\hspace{0.5in}&\text{Let}\; S\equiv\pm\{1, 5, 6, 7, 8, 9, 10, 11, 13, 14, 15, 19\}\pmod {42},\; \text{and} \\
&T\equiv\pm\{1, 4, 5, 6, 7, 9, 13, 14, 15, 17, 19, 20\} \pmod{42}.\; \text{Then}\\
&p(S,n)=p(T,n-1),\; \text{for all}\; n\geq 1.
\end{align*}
\begin{align*}
\mathrm{(ii)}\hspace{0.5in}&\text{Let}\; S\equiv\pm\{2, 3, 5, 7, 8, 9, 11, 12, 13, 14, 17, 19\}\pmod {42},\; \text{and} \\
&T\equiv\pm\{1, 3, 5, 7, 9, 11, 12, 14, 16, 17, 19, 20\} \pmod{42}.\; \text{Then}\\
&p(S,n)=p(T,n-2),\; \text{for all}\; n\geq 2.
\end{align*}
\begin{align*}
\mathrm{(iii)}\hspace{0.5in}&\text{Let}\; S\equiv\pm\{1, 3, 4, 5, 7, 11, 13, 14, 15, 16, 17, 18\}\pmod {42},\; \text{and} \\
&T\equiv\pm\{1, 2, 3, 7, 10, 11, 13, 14, 15, 17, 18, 19\} \pmod{42}.\; \text{Then}\\
&p(S,n)=p(T,n-1),\; \text{for all}\; n\geq 1.
\end{align*}
\end{theorem}
\begin{proof}
As we mentioned above the identities (i)--(iii) follow from
\eqref{qpp} or equivalently from \eqref{four2} with the following
sets of parameters
\begin{equation*}
[1, 3, 2, 7, 9],\;[1, 4, 2, 6, 9],\;[1, 2, 4, 12, 18]
\end{equation*}
and with $q$ replaced by $q^{21}$ in each instance.
\end{proof}

\begin{theorem}\label{42.2}

\begin{align*}
\mathrm{(i)}\hspace{0.5in}&\text{Let}\; S\equiv\pm\{1, 4, 5, 7, 8, 9, 11, 12, 13, 15, 16, 19\}\pmod {42},\; \text{and} \\
&T\equiv\pm\{1, 3, 4, 7, 8, 11, 12, 13, 15, 17, 19, 20\} \pmod{42}.\; \text{Then}\\
&p(S,n)=p(T,n-1),\; \text{for all}\; n\geq 1.
\end{align*}
\begin{align*}
\mathrm{(ii)}\hspace{0.5in}&\text{Let}\; S\equiv\pm\{2, 3, 4, 5, 7, 9, 11, 13, 17, 18, 19, 20\}\pmod {42},\; \text{and} \\
&T\equiv\pm\{1, 2, 5, 7, 9, 11, 13, 15, 16, 18, 19, 20\} \pmod{42}.\; \text{Then}\\
&p(S,n)=p(T,n-2),\; \text{for all}\; n\geq 2.
\end{align*}
\begin{align*}
\mathrm{(iii)}\hspace{0.5in}&\text{Let}\; S\equiv\pm\{1, 2, 3, 4, 5, 6, 7, 9, 10, 11, 17, 19\}\pmod {42},\; \text{and} \\
&T\equiv\pm\{1, 2, 3, 4, 5, 6, 7, 8, 11, 13, 15, 17\} \pmod{42}.\; \text{Then}\\
&p(S,n)=p(T,n),\; \text{for all}\; n\neq 8.
\end{align*}
\begin{align*}
\mathrm{(iv)}\hspace{0.5in}&\text{Let}\; S\equiv\pm\{1, 3, 5, 7, 8, 10, 11, 12, 13, 15, 17, 20\}\pmod {42},\; \text{and} \\
&T\equiv\pm\{1, 2, 5, 7, 9, 10, 12, 13, 15, 17, 19, 20\} \pmod{42}.\; \text{Then}\\
&p(S,n)=p(T,n-1),\; \text{for all}\; n\geq 1.
\end{align*}
\begin{align*}
\mathrm{(v)}\hspace{0.5in}&\text{Let}\; S\equiv\pm\{3, 4, 5, 6, 7, 9, 10, 11, 13, 16, 17, 19\}\pmod {42},\; \text{and} \\
&T\equiv\pm\{1, 3, 6, 7, 10, 11, 13, 15, 16, 17, 19, 20\} \pmod{42}.\; \text{Then}\\
&p(S,n)=p(T,n-3),\; \text{for all}\; n\geq 3.
\end{align*}
\begin{align*}
\mathrm{(vi)}\hspace{0.5in}&\text{Let}\; S\equiv\pm\{1, 3, 5, 7, 8, 9, 10, 11, 16, 17, 18, 19\}\pmod {42},\; \text{and} \\
&T\equiv\pm\{1, 2, 5, 7, 8, 9, 13, 15, 16, 17, 18, 19\} \pmod{42}.\; \text{Then}\\
&p(S,n)=p(T,n-1),\; \text{for all}\; n\geq 1.
\end{align*}
\end{theorem}
\begin{proof}
We prove (i) which we can express as
\begin{align}
&[3,17,20:42]-q[5,9,16:42]\notag\\
&=[1, 3, 4, 5, 7, 8, 9, 11, 12, 13, 15, 16, 17, 19,
20:42].\label{dd1}
\end{align}
Similarly the first part of Theorem \ref{42.1} is equivalent to
\begin{align}
&[4,17,20:42]-q[8,10,11:42]\notag\\
&=[1, 4, 5, 6, 7, 8, 9, 10, 11, 13, 14, 15, 17, 19,
20:42].\label{dd2}
\end{align}
Comparing these two identities, we find that \eqref{dd1} is
equivalent to
\begin{align}
&[3,12,16:42]\big\{[4,17,20:42]-q[8,10,11:42]\big\}\notag\\
&=[6,10,14:42]\big\{[3,17,20:42]-q[5,9,16:42]\big\}.\label{dd3}
\end{align}
Rearranging the terms of \eqref{dd3}, we arrive at
\begin{align}
&q[10,16:42]\big\{[5,6,9,14:42]-[3,8,11,12:42]\big\}\notag\\
&=[3,17:42]\big\{[6,10,14,20:42]-[4,12,16,20:42]\big\}.\label{dd4}
\end{align}
Employing \eqref{four} twice with the choice of parameters
$[1,3,6,9,12],$\\
$[1,3,7,13,17]$ and with $q$ replaced by $q^{42}$ in both
instances, we deduce
\begin{equation*}
[5,6,9,14:42]-[3,8,11,12:42]=q^3[2,3,6,17:42]\;\text{and}
\end{equation*}
\begin{equation*}
[6,10,14,20:42]-[4,12,16,20:42]=q^4[2,6,10,16:42].
\end{equation*}
These last two equations readily imply \eqref{dd4} and so the
proof of (i) is complete.
\end{proof}
\begin{theorem}\label{42.3}

\begin{align*}
\mathrm{(i)}\hspace{0.5in}&\text{Let}\; S\equiv\pm\{1, 3, 4, 5, 7, 9, 11, 16, 17, 18, 19, 20\}\pmod {42},\; \text{and} \\
&T\equiv\pm\{1, 2, 3, 7, 8, 11, 13, 15, 17, 18, 19, 20\} \pmod{42}.\; \text{Then}\\
&p(S,n)=p(T,n-1),\; \text{for all}\; n\geq 1.
\end{align*}
\begin{align*}
\mathrm{(ii)}\hspace{0.5in}&\text{Let}\; S\equiv\pm\{1, 2, 5, 6, 7, 9, 10, 11, 13, 15, 16, 19\}\pmod {42},\; \text{and} \\
&T\equiv\pm\{1, 3, 4, 5, 6, 7, 11, 13, 15, 16, 17, 20\} \pmod{42}.\; \text{Then}\\
&p(S,n)=p(T,n),\; \text{for all}\; n\neq 2.
\end{align*}
\begin{align*}
\mathrm{(iii)}\hspace{0.5in}&\text{Let}\; S\equiv\pm\{1, 3, 4, 5, 7, 9, 10, 11, 12, 17, 19, 20\}\pmod {42},\; \text{and} \\
&T\equiv\pm\{1, 2, 5, 7, 8, 9, 10, 11, 12, 13, 15, 19\} \pmod{42}.\; \text{Then}\\
&p(S,n)=p(T,n),\; \text{for all}\; n\neq 2.
\end{align*}
\begin{align*}
\mathrm{(iv)}\hspace{0.5in}&\text{Let}\; S\equiv\pm\{2, 3, 5, 7, 8, 9, 10, 11, 13, 17, 18, 19\}\pmod {42},\; \text{and} \\
&T\equiv\pm\{1, 3, 5, 7, 8, 11, 13, 15, 16, 17, 18, 20\} \pmod{42}.\; \text{Then}\\
&p(S,n)=p(T,n-2),\; \text{for all}\; n\geq 2.
\end{align*}
\begin{align*}
\mathrm{(v)}\hspace{0.5in}&\text{Let}\; S\equiv\pm\{3, 4, 5, 7, 8, 9, 10, 11, 12, 13, 17, 19\}\pmod {42},\; \text{and} \\
&T\equiv\pm\{1, 4, 5, 7, 9, 12, 13, 15, 16, 17, 19, 20\} \pmod{42}.\; \text{Then}\\
&p(S,n)=p(T,n-3),\; \text{for all}\; n\geq 3.
\end{align*}
\begin{align*}
\mathrm{(vi)}\hspace{0.5in}&\text{Let}\; S\equiv\pm\{1, 2, 4, 5, 6, 7, 9, 13, 15, 16, 17, 19\}\pmod {42},\; \text{and} \\
&T\equiv\pm\{1, 2, 3, 6, 7, 8, 10, 11, 13, 15, 17, 19\} \pmod{42}.\; \text{Then}\\
&p(S,n)=p(T,n),\; \text{for all}\; n\neq 3.
\end{align*}
\end{theorem}

\begin{proof}
We prove (i) which is equivalent to
\begin{align}
&[2,8,13,15:42]-q[4,5,9,16:42]\notag\\
&=[1, 2, 3, 4, 5, 7, 8, 9, 11, 13, 15, 16, 17, 18, 19,
20:42].\label{rdd1}
\end{align}
From  Theorem \ref{42.1}(iii), we also have
\begin{align}
&[2,10,19:42]-q[4,5,16:42]\notag\\
&=[1, 2, 3, 4, 5, 7, 10, 11, 13, 14, 15, 16, 17, 18,
19:42].\label{rdd2}
\end{align}
Therefore, it suffices to show that
\begin{align}
&[8,9,20:42]\big\{[2,10,19:42]-q[4,5,16:42]\big\}\notag\\
&=[10,14:42]\big\{[2,8,13,15:42]-q[4,5,9,16:42]\big\}.\label{rdd3}
\end{align}
After rearrangement of the terms of \eqref{dd3}, we arrive at
\begin{align}
&[2,8,10:42]\big\{[9,19,20:42]-[13,14,15:42]\big\}\notag\\
&=q[4,5,16,9:42]\big\{[8,20:42]-[10,14:42]\big\}.\label{rdd4}
\end{align}
Employing \eqref{four} twice with the choice of parameters
$[1,5,14,19,20],$\\
$[1,3,11,17,19]$ and with $q$ replaced by $q^{42}$ in both
instances, we find that
\begin{equation*}
[9,18,19,20:42]-[13,14,15,18:42]=-q^9[4,5,6,9:42]\;\text{and}
\end{equation*}
\begin{equation*}
[8,16,18,20:42]-[10,14,16,18:42]=-q^8[2,6,8,10:42].
\end{equation*}
Now, \eqref{rdd4} follows easily from these last two equations and
so the proof of (i) is complete.
\end{proof}

\section{Modulus $M = 46$}
There is only one  class  of identities under multiplication
by the group $U(46)$. \\
\begin{theorem}\label{46.1}

\begin{align*}
\mathrm{(i)}\hspace{0.5in}&\text{Let}\; S\equiv\pm\{1, 3, 4, 5, 7, 13, 15, 16, 17, 18, 19, 22\}\pmod {46},\; \text{and} \\
&T\equiv\pm\{1, 2, 3, 7, 11, 12, 15, 16, 17, 19, 20, 21\}\pmod{46}.\; \text{Then}\\
&p(S,n)=p(T,n-1),\; \text{for all}\; n\geq 1.
\end{align*}
\begin{align*}
\mathrm{(ii)}\hspace{0.5in}&\text{Let}\; S\equiv\pm\{1, 2, 3, 5, 7, 8, 9, 11, 12, 15, 20, 21\}\pmod {46},\; \text{and} \\
&T\equiv\pm\{1, 2, 3, 5, 6, 9, 10, 11, 13, 14, 17, 21\}\pmod{46}.\; \text{Then}\\
&p(S,n)=p(T,n),\; \text{for all}\; n\neq 6.
\end{align*}
\begin{align*}
\mathrm{(iii)}\hspace{0.5in}&\text{Let}\; S\equiv\pm\{3, 5, 7, 8, 9, 10, 11, 12, 13, 14, 15, 17\}\pmod {46},\; \text{and} \\
&T\equiv\pm\{2, 3, 5, 7, 11, 12, 15, 17, 18, 19, 20, 21\}\pmod{46}.\; \text{Then}\\
&p(S,n)=p(T,n-3),\; \text{for all}\; n\geq 3.
\end{align*}
\begin{align*}
\mathrm{(iv)}\hspace{0.5in}&\text{Let}\; S\equiv\pm\{2, 3, 5, 7, 8, 9, 13, 14, 15, 19, 20, 21\}\pmod {46},\; \text{and} \\
&T\equiv\pm\{1, 3, 5, 7, 11, 12, 13, 16, 18, 19, 20, 21\}\pmod{46}.\; \text{Then}\\
&p(S,n)=p(T,n-2),\; \text{for all}\; n\geq 2.
\end{align*}
\begin{align*}
\mathrm{(v)}\hspace{0.5in}&\text{Let}\; S\equiv\pm\{3, 4, 5, 6, 7, 9, 13, 15, 16, 17, 18, 19\}\pmod {46},\; \text{and} \\
&T\equiv\pm\{1, 3, 6, 9, 10, 13, 14, 15, 17, 19, 21, 22\}\pmod{46}.\; \text{Then}\\
&p(S,n)=p(T,n-3),\; \text{for all}\; n\geq 3.
\end{align*}
\begin{align*}
\mathrm{(vi)}\hspace{0.5in}&\text{Let}\; S\equiv\pm\{2, 3, 5, 8, 9, 11, 12, 13, 14, 15, 19, 21\}\pmod {46},\; \text{and} \\
&T\equiv\pm\{1, 3, 6, 8, 10, 11, 13, 15, 17, 19, 21, 22\}\pmod{46}.\; \text{Then}\\
&p(S,n)=p(T,n-2),\; \text{for all}\; n\geq 2.
\end{align*}
\begin{align*}
\mathrm{(vii)}\hspace{0.5in}&\text{Let}\; S\equiv\pm\{1, 4, 6, 7, 9, 10, 11, 13, 15, 17, 19, 22\}\pmod {46},\; \text{and} \\
&T\equiv\pm\{1, 3, 5, 7, 9, 11, 13, 16, 17, 18, 20, 22\}\pmod{46}.\; \text{Then}\\
&p(S,n)=p(T,n-1),\; \text{for all}\; n\geq 1.
\end{align*}
\begin{align*}
\mathrm{(viii)}\hspace{0.5in}&\text{Let}\; S\equiv\pm\{1, 5, 6, 8, 9, 10, 11, 13, 14, 15, 17, 21\}\pmod {46},\; \text{and} \\
&T\equiv\pm\{1, 4, 5, 7, 9, 10, 13, 15, 16, 19, 21, 22\}\pmod{46}.\; \text{Then}\\
&p(S,n)=p(T,n-1),\; \text{for all}\; n\geq 1.
\end{align*}
\begin{align*}
\mathrm{(ix)}\hspace{0.5in}&\text{Let}\; S\equiv\pm\{1, 4, 5, 6, 7, 9, 13, 16, 17, 19, 21, 22\}\pmod {46},\; \text{and} \\
&T\equiv\pm\{1, 3, 4, 5, 11, 12, 13, 17, 18, 19, 20, 21\}\pmod{46}.\; \text{Then}\\
&p(S,n)=p(T,n-1),\; \text{for all}\; n\geq 1.
\end{align*}
\begin{align*}
\mathrm{(x)}\hspace{0.5in}&\text{Let}\; S\equiv\pm\{1, 3, 4, 5, 7, 9, 11, 16, 17, 18, 19, 20\}\pmod {46},\; \text{and} \\
&T\equiv\pm\{1, 2, 5, 7, 8, 9, 11, 12, 15, 18, 19, 21\}\pmod{46}.\; \text{Then}\\
&p(S,n)=p(T,n),\; \text{for all}\; n\neq 2.
\end{align*}
\begin{align*}
\mathrm{(xi)}\hspace{0.5in}&\text{Let}\; S\equiv\pm\{2, 3, 7, 8, 9, 10, 11, 13, 14, 15, 17, 21\}\pmod {46},\; \text{and} \\
&T\equiv\pm\{1, 4, 6, 7, 9, 11, 14, 15, 17, 19, 21, 22\}\pmod{46}.\; \text{Then}\\
&p(S,n)=p(T,n-2),\; \text{for all}\; n\geq 2.
\end{align*}

\end{theorem}
\begin{proof}
The identities (i)--(xi) follow from \eqref{four2} with the
following choice of parameters
\begin{align*}
&[1, 2, 4, 12, 20],\;[1, 2, 8, 10, 11],\;[1, 6, 3, 9, 10],\; [1, 4,
2, 6, 9],\\&[1, 4, 5, 14, 20],\;[1, 4, 2, 6, 10],\;[1, 2, 5, 17,
18],\;[1, 3, 2, 7, 10],\\ &[1, 3, 2, 6, 8],[1, 2, 4, 8, 9], [1, 3,
4, 16, 18],
\end{align*}
and with $q$ replaced by $q^{23}$ in each instance.
\end{proof}

\section{Modulus $M = 48$}
There are seven distinct classes  of identities under multiplication
by the group $U(48)$. The first identity in Theorem \ref{48.1} below
was first given  by Kalvade \cite{KA}
and the second identity in Theorem \ref{48.3} was first given by Alladi \cite{AG}.\\
\begin{theorem}\label{48.1}

\begin{align*}
\mathrm{(i)}\hspace{0.5in}&\text{Let}\; S\equiv\pm\{1, 3, 4, 5, 7, 14, 16, 17, 18, 19, 21, 23\} \pmod {48},\; \text{and} \\
&T\equiv\pm\{1, 2, 3, 7, 11, 13, 16, 17, 18, 20, 21, 23\} \pmod{48}.\; \text{Then}\\
&p(S,n)=p(T,n-1),\; \text{for all}\; n\geq 1.
\end{align*}
\begin{align*}
\mathrm{(ii)}\hspace{0.5in}&\text{Let}\; S\equiv\pm\{4, 5, 6, 7, 9, 10, 11, 13, 15, 16, 17, 19\} \pmod {48},\; \text{and} \\
&T\equiv\pm\{1, 5, 6, 9, 11, 13, 15, 16, 19, 20, 22, 23\} \pmod{48}.\; \text{Then}\\
&p(S,n)=p(T,n-4),\; \text{for all}\; n\geq 4.
\end{align*}
\end{theorem}
\begin{proof}
The identities (i) and (ii) follow from \eqref{four2} with the
choice of parameters $[1,2,4,12,21]$,\\ $[1,6,2,8,10]$ and with $q$
replaced by $q^{24}$ in both instances.
\end{proof}

\begin{theorem}\label{48.2}

\begin{align*}
\mathrm{(i)}\hspace{0.5in}&\text{Let}\; S\equiv\pm\{1, 4, 5, 6, 7, 9, 15, 17, 19, 20, 22, 23\} \pmod {48},\; \text{and} \\
&T\equiv\pm\{1, 3, 4, 5, 11, 13, 14, 18, 19, 20, 21, 23\} \pmod{48}.\; \text{Then}\\
&p(S,n)=p(T,n-1),\; \text{for all}\; n\geq 1.
\end{align*}
\begin{align*}
\mathrm{(ii)}\hspace{0.5in}&\text{Let}\; S\equiv\pm\{2, 3, 4, 5, 7, 11, 13, 17, 18, 19, 20, 21\} \pmod {48},\; \text{and} \\
&T\equiv\pm\{1, 4, 6, 7, 9, 10, 11, 13, 15, 17, 20, 23\} \pmod{48}.\; \text{Then}\\
&p(S,n)=p(T,n),\; \text{for all}\; n\neq 1.
\end{align*}
\end{theorem}
\begin{proof}
The identities (i) and (ii) follow from \eqref{four2} with the
choice of parameters $[1,3,2,6,8]$,\\ $[1, 3, 4, 8, 10]$ and with
$q$ replaced by $q^{24}$ in both instances.
\end{proof}
\begin{theorem}\label{48.3}

\begin{align*}
\mathrm{(i)}\hspace{0.5in}&\text{Let}\; S\equiv\pm\{1, 6, 7, 8, 9, 10, 11, 13, 15, 16, 17, 23\} \pmod {48},\; \text{and} \\
&T\equiv\pm\{1, 5,6, 7, 8, 9, 15, 16, 17,19,22, 23\} \pmod{48}.\; \text{Then}\\
&p(S,n)=p(T,n-1),\; \text{for all}\; n\geq 1.
\end{align*}
\begin{align*}
\mathrm{(ii)}\hspace{0.5in}&\text{Let}\; S\equiv\pm\{2, 3, 5, 7, 8, 11, 13, 16, 17, 18, 19, 21\} \pmod {48},\; \text{and} \\
&T\equiv\pm\{1, 3, 5, 8, 11, 13, 14, 16, 18, 19, 21, 23\} \pmod{48}.\; \text{Then}\\
&p(S,n)=p(T,n-2),\; \text{for all}\; n\geq 2.
\end{align*}
\end{theorem}
\begin{proof}
The identities (i) and (ii) follow from \eqref{four2} with the
choice of parameters $[1,3,2,8,10]$,\\ $[1, 3, 4, 14, 20]$ and with
$q$ replaced by $q^{24}$ in both instances.
\end{proof}
\begin{theorem}\label{48.4}

\begin{align*}
\hspace{0.5in}&\text{Let}\; S\equiv\pm\{3, 4, 5, 6, 7, 9, 15, 16, 17, 18, 19, 21\} \pmod {48},\; \text{and} \\
&T\equiv\pm\{1, 3, 6, 9, 11, 13, 15, 16, 18, 20, 21, 23\} \pmod{48}.\; \text{Then}\\
&p(S,n)=p(T,n-3),\; \text{for all}\; n\geq 3.
\end{align*}
\end{theorem}
\begin{proof}
This identity follows from  \eqref{four2} with the choice of
parameters $[1, 4, 5, 14, 21]$ and with $q$ replaced by $q^{24}$.
\end{proof}
\begin{theorem}\label{48.5}

\begin{align*}
\mathrm{(i)}\hspace{0.5in}&\text{Let}\; S\equiv\pm\{1, 5, 6, 8, 9, 10, 11, 13, 15, 19, 20, 23\} \pmod {48},\; \text{and} \\
&T\equiv\pm\{1, 4, 5, 7, 9, 10, 15, 17, 18, 19, 20, 23\} \pmod{48}.\; \text{Then}\\
&p(S,n)=p(T,n-1),\; \text{for all}\; n\geq 1.
\end{align*}
\begin{align*}
\mathrm{(ii)}\hspace{0.5in}&\text{Let}\; S\equiv\pm\{1, 2, 3, 4, 5, 6, 11, 13, 19, 20, 21, 23\} \pmod {48},\; \text{and} \\
&T\equiv\pm\{1, 2, 3, 4, 5, 7, 8, 17, 18, 19, 21, 23\} \pmod{48}.\; \text{Then}\\
&p(S,n)=p(T,n),\; \text{for all}\; n\neq 6.
\end{align*}
\begin{align*}
\mathrm{(iii)}\hspace{0.5in}&\text{Let}\; S\equiv\pm\{4, 5, 6, 7, 8, 9, 11, 13, 15, 17, 19, 22\} \pmod {48},\; \text{and} \\
&T\equiv\pm\{1, 4, 7, 9, 11, 13, 15, 17, 18, 20, 22, 23\} \pmod{48}.\; \text{Then}\\
&p(S,n)=p(T,n-4),\; \text{for all}\; n\geq 4.
\end{align*}
\begin{align*}
\mathrm{(iv)}\hspace{0.5in}&\text{Let}\; S\equiv\pm\{3, 4, 5, 6, 7, 11, 13, 14, 17, 19, 20, 21\} \pmod {48},\; \text{and} \\
&T\equiv\pm\{1, 3, 7, 8, 11, 13, 14, 17, 18, 20, 21, 23\} \pmod{48}.\; \text{Then}\\
&p(S,n)=p(T,n-3),\; \text{for all}\; n\geq 3.
\end{align*}
\end{theorem}
\begin{proof}
We prove (i) which is  equivalent to
\begin{align*}
&[4,7,17,18:48]-q[6,8,11,13:48]\\&=[1, 4, 5, 6, 7, 8, 9, 10, 11,
13, 15, 17, 18, 19, 20, 23:48].
\end{align*}
From Theorem \ref{48.1}(ii), we similarly have that
\begin{align*}
&[1,20,22,23:48]-q^4[4,7,10,17:48]\\& =[1, 4, 5, 6, 7, 9, 10, 11,
13, 15, 16, 17, 19, 20, 22, 23:48].
\end{align*}
Therefore, it suffices to prove that
\begin{align}
&[8,18:48]\bigr\{[1,20,22,23:48]-q^4[4,7,10,17:48]\bigl\}\notag\\
&=[16,22:48]\bigr\{[4,7,17,18:48]-q[6,8,11,13:48]\bigl\}.\label{c2}
\end{align}
Rearranging the  terms in \eqref{c2}, we derive the equivalent
identity
\begin{align}
&[7,17:48]\bigr\{[4,16,22,18:48]+q^4[4,8,10,18:48]\bigl\}\notag\\
&=[8,22:48]\bigr\{[1,18,20,23:48]+q[6,11,13,16:48]\bigl\}.\label{c3}
\end{align}
We now employ \eqref{four} twice with the set of variables
$[1,5,9,17,19]$, $[1,7,8,19,21]$ and with $q$ replaced by $q^{48}$
in each instance, we find that
\begin{equation}\label{c4}
[4,16,22,18:48]+q^4[4,8,10,18:48]=[8,12,14,22:48]\; \text{and}
\end{equation}
\begin{equation}\label{c5}
[1,18,20,23:48]+q[6,11,13,16:48]=[7,12,14,17:48].
\end{equation}
We see that \eqref{c4} together with \eqref{c5} clearly implies
\eqref{c3} and so the proof of (i) is complete.
\end{proof}
\begin{theorem}\label{48.6}

\begin{align*}
\mathrm{(i)}\hspace{0.5in}&\text{Let}\; S\equiv\pm\{1, 4, 7, 9, 10, 11, 12, 13, 14, 15, 17, 23\} \pmod {48},\; \text{and} \\
&T\equiv\pm\{1, 3, 7, 8, 10, 11, 12, 13, 17, 21, 22, 23\} \pmod{48}.\; \text{Then}\\
&p(S,n)=p(T,n-1),\; \text{for all}\; n\geq 1.
\end{align*}
\begin{align*}
\mathrm{(ii)}\hspace{0.5in}&\text{Let}\; S\equiv\pm\{2, 5, 7, 8, 9, 11, 12, 13, 14, 15, 17, 19\} \pmod {48},\; \text{and} \\
&T\equiv\pm\{2, 3, 5, 7, 11, 12, 13, 17, 19, 20, 21, 22\} \pmod{48}.\; \text{Then}\\
&p(S,n)=p(T,n-2),\; \text{for all}\; n\geq 2.
\end{align*}
\begin{align*}
\mathrm{(iii)}\hspace{0.5in}&\text{Let}\; S\equiv\pm\{1, 3, 5, 7, 8, 10, 12, 17, 19, 21, 22, 23\} \pmod {48},\; \text{and} \\
&T\equiv\pm\{1, 2, 5, 7, 9, 12, 15, 17, 19, 20, 22, 23\} \pmod{48}.\; \text{Then}\\
&p(S,n)=p(T,n-1),\; \text{for all}\; n\geq 1.
\end{align*}
\begin{align*}
\mathrm{(iv)}\hspace{0.5in}&\text{Let}\; S\equiv\pm\{1, 2, 5, 8, 9, 11, 12, 13, 14, 15, 19, 23\} \pmod {48},\; \text{and} \\
&T\equiv\pm\{1, 3, 4, 5, 10, 11, 12, 13, 14, 19, 21, 23\} \pmod{48}.\; \text{Then}\\
&p(S,n)=p(T,n),\; \text{for all}\; n\neq 2.
\end{align*}
\end{theorem}
\begin{proof}
We prove (i). First we write (i) in its equivalent form
\begin{align}
&[3,8,21,22:48]-q[4,9,14,15:48]\notag\\
&=[1, 3, 4, 7, 8, 9, 10, 11, 12, 13, 14, 15, 17, 21, 22,
23:48].\label{cc1}
\end{align}
Similarly from the first part of Theorem \ref{48.1}
\begin{align}
&[2,11,13,20:48]-q[4,5,14,19:48]\notag\\
&=[1, 2, 3, 4, 5, 7, 11, 13, 14, 16, 17, 18, 19, 20, 21,
23:48].\label{cc2}
\end{align}
Therefore, it suffices to show that
\begin{align}
&[8,9,10,12,15,22:48]\big\{[2,11,13,20:48]-q[4,5,14,19:48]\big\}\notag\\
&=[2,5,16,18,19,20:48]\big\{[3,8,21,22:48]-q[4,9,14,15:48]\big\}.\label{cc3}
\end{align}
After rearrangement of the terms, we arrive at
\begin{align}
&[2,22,8,20:48]\big\{[9,10,11,12,13,15:48]-[3,5,16,18,19,21:48]\big\}\notag\\
&=-q[4,5,9,14,15,19:48]\big\{[2,16,18,20:48]-[8,10,12,22:48]\big\}.\label{cc4}
\end{align}
By \eqref{four} with the choice of parameters $[1,7,9,17,19]$ and
with $q$ replaced by $q^{48}$, we deduce that
\begin{equation}\label{cc5}
[2,16,18,20:48]-[8,10,12,22:48]=-q^2[6,8,10,20:48].
\end{equation}
Using \eqref{cc5} in \eqref{cc4}, we arrive at
\begin{align}
&[2,22:48]\big\{[3,5,16,18,19,21:48]-[9,10,11,12,13,15:48]\big\}\notag\\
&=-q^3[4,5,6,9,10,14,15,19:48].\label{cc6}
\end{align}
Switching to base $24$, we have from \eqref{cc6} that
\begin{align}
&[2:24]\big\{[3,5,8,9:24](8,9:24)-[5,6,9,11:24](5,6:24)\big\}\notag\\
&=-q^3[2,3,5,9,10:24](2,3:24),\label{cc7}
\end{align}
which simplifies to
\begin{align}
&[3,8:24](8,9:24)-[6,11:24](5,6:24)\notag\\
&=-q^3[3,10:24](2,3:24).\label{cc8}
\end{align}
By employing \eqref{four} with $a,\,b,\,c,\,x,\,y$ and $q$
replaced by $-q,\,q^3,\,q^6,\,-q^9,\,q^9$ and $q^{24}$,
respectively, we establish \eqref{cc8} and complete the proof of
(i).
\end{proof}
\begin{theorem}\label{48.7}

\begin{align*}
\mathrm{(i)}\hspace{0.5in}&\text{Let}\; S\equiv\pm\{1, 4, 5, 7, 8, 11, 14, 15, 17, 18, 19, 21\} \pmod {48},\; \text{and} \\
&T\equiv\pm\{1, 3, 4, 7, 11, 13, 14, 15, 18, 19, 20, 23\}\pmod{48}.\; \text{Then}\\
&p(S,n)=p(T,n-1),\; \text{for all}\; n\geq 1.
\end{align*}
\begin{align*}
\mathrm{(ii)}\hspace{0.5in}&\text{Let}\; S\equiv\pm\{1, 5, 6, 7, 8, 9, 11, 13, 20, 21, 22, 23\} \pmod {48},\; \text{and} \\
&T\equiv\pm\{1, 4, 5, 6, 7, 13, 15, 17, 19, 20, 21, 22\}\pmod{48}.\; \text{Then}\\
&p(S,n)=p(T,n-1),\; \text{for all}\; n\geq 1.
\end{align*}
\begin{align*}
\mathrm{(iii)}\hspace{0.5in}&\text{Let}\; S\equiv\pm\{1, 2, 4, 5, 7, 9, 11, 17, 18, 19, 20, 21\} \pmod {48},\; \text{and} \\
&T\equiv\pm\{1, 2, 3, 7, 8, 9, 11, 13, 18, 19, 20, 23\}\pmod{48}.\; \text{Then}\\
&p(S,n)=p(T,n),\; \text{for all}\; n\neq 3.
\end{align*}
\begin{align*}
\mathrm{(iv)}\hspace{0.5in}&\text{Let}\; S\equiv\pm\{4, 5, 6, 7, 8, 9, 10, 11, 17, 19, 21, 23\} \pmod {48},\; \text{and} \\
&T\equiv\pm\{1, 4, 6, 10, 11, 13, 15, 17, 19, 20, 21, 23\}\pmod{48}.\; \text{Then}\\
&p(S,n)=p(T,n-4),\; \text{for all}\; n\geq 4.
\end{align*}
\begin{align*}
\mathrm{(v)}\hspace{0.5in}&\text{Let}\; S\equiv\pm\{1, 3, 4, 5, 6, 7, 9, 10, 11, 13, 20, 23\} \pmod {48},\; \text{and} \\
&T\equiv\pm\{1, 3, 4, 5, 6, 7, 8, 10, 13, 15, 17, 19\}\pmod{48}.\; \text{Then}\\
&p(S,n)=p(T,n),\; \text{for all}\; n\neq 8.
\end{align*}
\begin{align*}
\mathrm{(vi)}\hspace{0.5in}&\text{Let}\; S\equiv\pm\{2, 3, 4, 5, 7, 13, 15, 17, 18, 19, 20, 23\} \pmod {48},\; \text{and} \\
&T\equiv\pm\{1, 2, 5, 8, 11, 13, 15, 17, 18, 20, 21, 23\}\pmod{48}.\; \text{Then}\\
&p(S,n)=p(T,n-2),\; \text{for all}\; n\geq 2.
\end{align*}
\begin{align*}
\mathrm{(vii)}\hspace{0.5in}&\text{Let}\; S\equiv\pm\{3, 4, 5, 6, 7, 9, 11, 17, 19, 20, 22, 23\} \pmod {48},\; \text{and} \\
&T\equiv\pm\{1, 3, 6, 8, 11, 13, 15, 17, 19, 20, 22, 23\}\pmod{48}.\; \text{Then}\\
&p(S,n)=p(T,n-3),\; \text{for all}\; n\geq 3.
\end{align*}
\begin{align*}
\mathrm{(viii)}\hspace{0.5in}&\text{Let}\; S\equiv\pm\{3, 4, 5, 7, 8, 9, 13, 14, 17, 18, 19, 23\} \pmod {48},\; \text{and} \\
&T\equiv\pm\{1, 4, 5, 9, 11, 13, 14, 17, 18, 20, 21, 23\}\pmod{48}.\; \text{Then}\\
&p(S,n)=p(T,n-3),\; \text{for all}\; n\geq 3.
\end{align*}
\end{theorem}
\begin{proof}
We prove (i). Recall that by \eqref{cc2}, we have

\begin{align}
&[2,11,13,20:48]-q[4,5,14,19:48]\notag\\
&=[1, 2, 3, 4, 5, 7, 11, 13, 14, 16, 17, 18, 19, 20, 21,
23:48].\label{ccc2}
\end{align}
We write (i) in its equivalent form
\begin{align}
&[3,13,20,23:48]-q[5,8,17,21:48]\notag\\
&=[1, 3, 4, 5, 7, 8, 11, 13, 14, 15, 17, 18, 19, 20, 21,
23:48].\label{ccc3}
\end{align}
Therefore, it suffices to prove that
\begin{align}
&[8,15:48]\big\{[2,11,13,20:48]-q[4,5,14,19:48]\big\}\notag\\
&=[2,16:48]\big\{[3,13,20,23:48]-q[5,8,17,21:48]\big\}.\label{ccc4}
\end{align}
By rearranging the terms in \eqref{ccc4}, we arrive at
\begin{align}
&[2,13:48]\big\{[8,11,15,20:48]-[3,16,20,23:48]\big\}\notag\\
&=q[5,8:48]\big\{[4,14,15,19:48]-[2,16,17,21:48]\big\}.\label{ccc5}
\end{align}
Next, we employ \eqref{four} twice with the choice of parameters
$[1,6,9,17,21]$,\\
$[1,3,5,17,18]$ and with $q$ replaced by $q^{48}$, we find that
\begin{equation*}
[8,11,15,20:48]-[3,16,20,23:48]=q^3[5,8,12,17:48] \;\text{and}
\end{equation*}
\begin{equation*}
[4,14,15,19:48]-[2,16,17,21:48]=q^2[2,12,13,17:48].
\end{equation*}
From these last two equations, \eqref{ccc5} readily follows and so
the proof of (i) is complete.
\end{proof}

\section{Modulus $M = 50$}
There are four distinct classes  of identities under multiplication
by the group $U(50)$. \\
\begin{theorem}\label{50.1}

\begin{align*}
\mathrm{(i)}\hspace{0.5in}&\text{Let}\; S\equiv\pm\{1, 3, 5, 7, 8, 13,15, 17, 18, 19, 20, 22\} \pmod {50},\; \text{and} \\
&T\equiv\pm\{1, 2, 5, 7, 12, 13, 15, 17, 18, 20, 21, 23\} \pmod{50}.\; \text{Then}\\
&p(S,n)=p(T,n-1),\; \text{for all}\; n\geq 1.
\end{align*}
\begin{align*}
\mathrm{(ii)}\hspace{0.5in}&\text{Let}\; S\equiv\pm\{1, 3, 4, 5, 6, 10, 11, 13, 14, 15, 19, 21\} \pmod {50},\; \text{and} \\
&T\equiv\pm\{1, 3, 4, 5, 7, 9, 10, 11, 15, 16, 21, 24\} \pmod{50}.\; \text{Then}\\
&p(S,n)=p(T,n),\; \text{for all}\; n\neq 6.
\end{align*}
\begin{align*}
\mathrm{(iii)}\hspace{0.5in}&\text{Let}\; S\equiv\pm\{1, 3, 5, 7, 9, 10, 11, 14, 15, 16, 19, 24\} \pmod {50},\; \text{and} \\
&T\equiv\pm\{1, 4, 5, 6, 7, 9, 10, 15, 17, 19, 21, 24\}\pmod{50}.\; \text{Then}\\
&p(S,n)=p(T,n),\; \text{for all}\; n\neq 3.
\end{align*}
\begin{align*}
\mathrm{(iv)}\hspace{0.5in}&\text{Let}\; S\equiv\pm\{3, 5, 7, 8, 9, 11, 12, 13, 15, 17, 18, 20\} \pmod {50},\; \text{and} \\
&T\equiv\pm\{2, 3, 5, 9, 12, 13, 15, 17, 20, 21, 22, 23\}\pmod{50}.\; \text{Then}\\
&p(S,n)=p(T,n-3),\; \text{for all}\; n\geq 3.
\end{align*}
\begin{align*}
\mathrm{(v)}\hspace{0.5in}&\text{Let}\; S\equiv\pm\{2, 5, 7, 8, 9, 11, 12, 13, 15, 17, 20, 23\} \pmod {50},\; \text{and} \\
&T\equiv\pm\{2, 3, 5, 7, 11, 13, 15, 18, 19, 20, 22, 23\}\pmod{50}.\; \text{Then}\\
&p(S,n)=p(T,n-2),\; \text{for all}\; n\geq 2.
\end{align*}
\begin{align*}
\mathrm{(vi)}\hspace{0.5in}&\text{Let}\; S\equiv\pm\{3, 4, 5, 9, 10, 11, 13, 14, 15, 16, 19, 21\} \pmod {50},\; \text{and} \\
&T\equiv\pm\{1, 5, 6, 9, 10, 13, 15, 16, 19, 21, 23, 24\}\pmod{50}.\; \text{Then}\\
&p(S,n)=p(T,n-3),\; \text{for all}\; n\geq 3.
\end{align*}
\begin{align*}
\mathrm{(vii)}\hspace{0.5in}&\text{Let}\; S\equiv\pm\{4, 5, 6, 7, 9, 10, 11, 15, 16, 17, 19, 21\} \pmod {50},\; \text{and} \\
&T\equiv\pm\{1, 5, 6, 10, 11, 14, 15, 17, 19, 21, 23, 24\}\pmod{50}.\; \text{Then}\\
&p(S,n)=p(T,n-4),\; \text{for all}\; n\geq 4.
\end{align*}
\begin{align*}
\mathrm{(viii)}\hspace{0.5in}&\text{Let}\; S\equiv\pm\{2, 3, 5, 7, 8, 11, 15, 17, 18, 19, 20, 23\} \pmod {50},\; \text{and} \\
&T\equiv\pm\{1, 3, 5, 8, 12, 13, 15, 17, 19, 20, 22, 23\}\pmod{50}.\; \text{Then}\\
&p(S,n)=p(T,n-2),\; \text{for all}\; n\geq 2.
\end{align*}
\begin{align*}
\mathrm{(ix)}\hspace{0.5in}&\text{Let}\; S\equiv\pm\{2, 3, 5, 7, 8, 9, 15, 17, 20, 21, 22, 23\} \pmod {50},\; \text{and} \\
&T\equiv\pm\{1, 3, 5, 7, 12, 13, 15, 18, 20, 21, 22, 23\}\pmod{50}.\; \text{Then}\\
&p(S,n)=p(T,n-2),\; \text{for all}\; n\geq 2.
\end{align*}
\begin{align*}
\mathrm{(x)}\hspace{0.5in}&\text{Let}\; S\equiv\pm\{1, 5, 6, 9, 10, 11, 13, 14, 15, 16, 19, 23\} \pmod {50},\; \text{and} \\
&T\equiv\pm\{1, 4, 5, 9, 10, 11, 14, 15, 17, 21, 23, 24\}\pmod{50}.\; \text{Then}\\
&p(S,n)=p(T,n-1),\; \text{for all}\; n\geq 1.
\end{align*}
\end{theorem}
\begin{proof}
To prove the identities  (i)--(x), we employ \eqref{four2}
with $q$ replaced by $q^{25}$ and with the following choice of
parameters for $a,\,b,\,c,\,x,$ and $y$ in each instance

$[1, 2, 4, 14, 21],[1, 2, 8, 11, 12],[1, 2, 5, 10, 11],[1, 4, 6,
18, 21],[1, 3, 6, 19, 20]$,\\$[1, 5, 2, 7, 11], [1, 6, 2, 8,
10],[1, 3, 4, 14, 21],[1, 4, 2, 6, 9]$, and $[1, 3, 2, 7, 11]$.

\end{proof}
\begin{theorem}\label{50.2}

\begin{align*}
\mathrm{(i)}\hspace{0.5in}&\text{Let}\; S\equiv\pm\{1, 6, 7, 8, 9, 11, 12, 13, 16, 17, 19, 23\} \pmod {50},\; \text{and} \\
&T\equiv\pm\{1, 5, 6, 7, 8, 11, 15, 17, 18, 19, 23, 24\} \pmod{50}.\; \text{Then}\\
&p(S,n)=p(T,n-1),\; \text{for all}\; n\geq 1.
\end{align*}
\begin{align*}
\mathrm{(ii)}\hspace{0.5in}&\text{Let}\; S\equiv\pm\{1, 3, 4, 5, 7, 15, 17, 18, 19, 21, 22, 24\} \pmod {50},\; \text{and} \\
&T\equiv\pm\{1, 2, 3, 7, 11, 14, 17, 18, 19, 21, 23, 24\} \pmod{50}.\; \text{Then}\\
&p(S,n)=p(T,n-1),\; \text{for all}\; n\geq 1.
\end{align*}
\begin{align*}
\mathrm{(iii)}\hspace{0.5in}&\text{Let}\; S\equiv\pm\{1, 3, 4, 5, 7, 9, 12, 13, 15, 16, 21, 22\} \pmod {50},\; \text{and} \\
&T\equiv\pm\{1, 3, 4, 6, 7, 8, 9, 13, 17, 19, 21, 22\} \pmod{50}.\; \text{Then}\\
&p(S,n)=p(T,n),\; \text{for all}\; n\neq 5.
\end{align*}
\begin{align*}
\mathrm{(iv)}\hspace{0.5in}&\text{Let}\; S\equiv\pm\{2, 3, 5, 9, 11, 12, 13, 14, 15, 16, 21, 23\} \pmod {50},\; \text{and} \\
&T\equiv\pm\{1, 3, 7, 9, 11, 12, 13, 16, 18, 21, 23, 24\} \pmod{50}.\; \text{Then}\\
&p(S,n)=p(T,n-2),\; \text{for all}\; n\geq 2.
\end{align*}
\begin{align*}
\mathrm{(v)}\hspace{0.5in}&\text{Let}\; S\equiv\pm\{2, 5, 6, 8, 9, 11, 13, 14, 15, 17, 19, 23\} \pmod {50},\; \text{and} \\
&T\equiv\pm\{2, 3, 4, 9, 11, 13, 14, 17, 19, 21, 22, 23\} \pmod{50}.\; \text{Then}\\
&p(S,n)=p(T,n-2),\; \text{for all}\; n\geq 2.
\end{align*}
\end{theorem}
\begin{proof}
Identities  (i)--(x) also follow from \eqref{four2}
 with $q$ replaced by $q^{25}$ and with the following
choice of parameters for $a,\,b,\,c,\,x,$ and $y$ in each
instance

$[1, 3, 2, 8, 10],[1, 2, 4, 12, 22],[1, 2, 7, 10, 11],[1, 4, 2, 6,
11]$, and $[1, 3, 6, 18, 20].$
\end{proof}

\begin{theorem}\label{50.3}
\begin{align*}
\mathrm{(i)}\hspace{0.5in}&\text{Let}\; S\equiv\pm\{1, 3, 5, 7, 8, 12, 13, 17, 19, 22, 23, 24\} \pmod {50},\; \text{and} \\
&T\equiv\pm\{1, 2, 5, 7, 11, 12, 15, 18, 19, 21, 23, 24\} \pmod{50}.\; \text{Then}\\
&p(S,n)=p(T,n-1),\; \text{for all}\; n\geq 1.
\end{align*}
\begin{align*}
\mathrm{(ii)}\hspace{0.5in}&\text{Let}\; S\equiv\pm\{3, 4, 5, 6, 7, 13, 14, 15, 17, 19, 21, 22\} \pmod {50},\; \text{and} \\
&T\equiv\pm\{1, 3, 7, 9, 11, 14, 15, 16, 19, 21, 22, 24\} \pmod{50}.\; \text{Then}\\
&p(S,n)=p(T,n-3),\; \text{for all}\; n\geq 3.
\end{align*}
\begin{align*}
\mathrm{(iii)}\hspace{0.5in}&\text{Let}\; S\equiv\pm\{1, 4, 6, 7, 9, 11, 15, 16, 17, 18, 19, 21\} \pmod {50},\; \text{and} \\
&T\equiv\pm\{1, 3, 5, 7, 11, 14, 15, 16, 17, 18, 23, 24\} \pmod{50}.\; \text{Then}\\
&p(S,n)=p(T,n-1),\; \text{for all}\; n\geq 1.
\end{align*}
\begin{align*}
\mathrm{(iv)}\hspace{0.5in}&\text{Let}\; S\equiv\pm\{2, 3, 5, 7, 8, 9, 13, 16, 17, 21, 22, 23\} \pmod {50},\; \text{and} \\
&T\equiv\pm\{1, 5, 7, 8, 9, 11, 12, 13, 15, 16, 18, 21\} \pmod{50}.\; \text{Then}\\
&p(S,n)=p(T,n),\; \text{for all}\; n\neq 1.
\end{align*}
\begin{align*}
\mathrm{(v)}\hspace{0.5in}&\text{Let}\; S\equiv\pm\{3, 5, 7, 8, 9, 11, 12, 13, 14, 17, 18, 23\} \pmod {50},\; \text{and} \\
&T\equiv\pm\{2, 3, 5, 9, 11, 14, 15, 18, 19, 21, 22, 23\} \pmod{50}.\; \text{Then}\\
&p(S,n)=p(T,n-3),\; \text{for all}\; n\geq 3.
\end{align*}
\begin{align*}
\mathrm{(vi)}\hspace{0.5in}&\text{Let}\; S\equiv\pm\{1, 3, 4, 6, 9, 11, 12, 13, 14, 15, 19, 21\} \pmod {50},\; \text{and} \\
&T\equiv\pm\{1, 3, 5, 6, 7, 9, 12, 13, 15, 16, 23, 24\} \pmod{50}.\; \text{Then}\\
&p(S,n)=p(T,n),\; \text{for all}\; n\neq 4.
\end{align*}
\begin{align*}
\mathrm{(vii)}\hspace{0.5in}&\text{Let}\; S\equiv\pm\{4, 5, 6, 7, 8, 9, 13, 15, 16, 17, 19, 23\} \pmod {50},\; \text{and} \\
&T\equiv\pm\{1, 4, 8, 9, 11, 14, 15, 17, 19, 21, 23, 24\} \pmod{50}.\; \text{Then}\\
&p(S,n)=p(T,n-4),\; \text{for all}\; n\geq 4.
\end{align*}
\begin{align*}
\mathrm{(viii)}\hspace{0.5in}&\text{Let}\; S\equiv\pm\{2, 3, 5, 6, 7, 11, 13, 17, 18, 19, 22, 23\} \pmod {50},\; \text{and} \\
&T\equiv\pm\{1, 5, 6, 8, 9, 11, 12, 13, 15, 17, 19, 22\} \pmod{50}.\; \text{Then}\\
&p(S,n)=p(T,n),\; \text{for all}\; n\neq 1.
\end{align*}
\begin{align*}
\mathrm{(xi)}\hspace{0.5in}&\text{Let}\; S\equiv\pm\{1, 2, 3, 4, 5, 7, 12, 13, 17, 18, 21, 23\} \pmod {50},\; \text{and} \\
&T\equiv\pm\{1, 2, 3, 4, 5, 8, 9, 15, 17, 19, 21, 22\}\pmod{50}.\; \text{Then}\\
&p(S,n)=p(T,n),\; \text{for all}\; n\neq 7.
\end{align*}
\begin{align*}
\mathrm{(x)}\hspace{0.5in}&\text{Let}\; S\equiv\pm\{2, 3, 4, 5, 11, 13, 14, 15, 17, 21, 23, 24\} \pmod {50},\; \text{and} \\
&T\equiv\pm\{1, 2, 6, 9, 11, 13, 15, 16, 19, 21, 23, 24\}\pmod{50}.\; \text{Then}\\
&p(S,n)=p(T,n-2),\; \text{for all}\; n\geq 2.
\end{align*}
\end{theorem}
\begin{proof}
We prove (i) which is  equivalent to
\begin{align*}
&[2,11,15,18,21:50]-q[3,8,13,17,22:50]\\&=[1,2,3,5,7,8,11,12,13,15,17,18,19,21,22,23,24:50].
\end{align*}
From  Theorem \ref{50.1}(i), we similarly have that
\begin{align*}
&[2,12,21,23:50]-q[3,8,19,22:50]\\&
=[1,2,3,5,7,8,12,13,15,17,18,19,20,21,22,23:50].
\end{align*}
Therefore, it suffices to prove that
\begin{align}
&[11,24:50]\bigr\{[2,12,21,23:50]-q[3,8,19,22:50]\bigl\}\notag\\
&=[20:50]\bigr\{[2,11,15,18,21:50]-q[3,8,13,17,22:50]\bigl\}.\label{st2}
\end{align}
After regrouping terms in \eqref{st2}, we are led to verify that
\begin{align}
&[21,2,11]\bigr\{[24,12,23:50]-[20,15,18:50]\bigl\}\notag\\
&=-q[3,8,22:50]\bigr\{[20,13,17:50]-[11,19,24:50]\bigl\}.\label{st5}
\end{align}
In \eqref{four}, we replace $q$ by $q^{50}$ and take $a=1$, $b=4$,
$c=16$, $x=22$, and $y=24$, we find that
\begin{equation}\label{st6}
[15,18,20,21:50]-[12,21,23,24:50]=q^{12}[3,6,8,9:50].
\end{equation}
Similarly, with the choice of parameters $a=1$, $b=3$, $c=14$,
$x=20$, and $y=23$, we deduce that
\begin{equation}\label{st7}
[13,17,20,22:50]-[11,19,22,24:50]=q^{11}[2,6,9,11:50].
\end{equation}
Equation \eqref{st6} together with \eqref{st7} clearly implies
\eqref{st5} and so the proof of (i) is complete.
\end{proof}
\begin{theorem}\label{50.4}
\begin{align*}
\mathrm{(i)}\hspace{0.5in}&\text{Let}\; S\equiv\pm\{1, 3, 5, 7, 11, 12, 13, 16, 17, 18, 22, 23\} \pmod {50},\; \text{and} \\
&T\equiv\pm\{1, 2, 5, 9, 11, 12, 15, 16, 17, 21, 22, 23\} \pmod{50}.\; \text{Then}\\
&p(S,n)=p(T,n-1),\; \text{for all}\; n\geq 1.
\end{align*}
\begin{align*}
\mathrm{(ii)}\hspace{0.5in}&\text{Let}\; S\equiv\pm\{1, 2, 3, 5, 6, 13, 14, 15, 16, 17, 19, 23\} \pmod {50},\; \text{and} \\
&T\equiv\pm\{1, 2, 3, 4, 9, 11, 14, 15, 16, 17, 19, 21\} \pmod{50}.\; \text{Then}\\
&p(S,n)=p(T,n),\; \text{for all}\; n\neq 4.
\end{align*}
\begin{align*}
\mathrm{(iii)}\hspace{0.5in}&\text{Let}\; S\equiv\pm\{3, 4, 5, 7, 11, 12, 13, 14, 15, 16, 19, 23\} \pmod {50},\; \text{and} \\
&T\equiv\pm\{1, 4, 7, 9, 11, 12, 15, 16, 19, 21, 23, 24\} \pmod{50}.\; \text{Then}\\
&p(S,n)=p(T,n-3),\; \text{for all}\; n\geq 3.
\end{align*}
\begin{align*}
\mathrm{(iv)}\hspace{0.5in}&\text{Let}\; S\equiv\pm\{1, 2, 3, 5, 6, 7, 8, 9, 11, 15, 18, 19\} \pmod {50},\; \text{and} \\
&T\equiv\pm\{1, 2, 3, 5, 6, 7, 8, 9, 12, 13, 17, 23\} \pmod{50}.\; \text{Then}\\
&p(S,n)=p(T,n),\; \text{for all}\; n\neq 11.
\end{align*}
\begin{align*}
\mathrm{(v)}\hspace{0.5in}&\text{Let}\; S\equiv\pm\{2, 3, 5, 7, 8, 11, 13, 17, 18, 21, 23, 24\} \pmod {50},\; \text{and} \\
&T\equiv\pm\{1, 3, 5, 8, 11, 13, 15, 18, 19, 21, 22, 24\} \pmod{50}.\; \text{Then}\\
&p(S,n)=p(T,n-2),\; \text{for all}\; n\geq 2.
\end{align*}
\begin{align*}
\mathrm{(vi)}\hspace{0.5in}&\text{Let}\; S\equiv\pm\{1, 6, 7, 8, 9, 11, 13, 14, 15, 16, 19, 21\} \pmod {50},\; \text{and} \\
&T\equiv\pm\{1, 5, 6, 7, 8, 13, 14, 15, 17, 21, 23, 24\} \pmod{50}.\; \text{Then}\\
&p(S,n)=p(T,n-1),\; \text{for all}\; n\geq 1.
\end{align*}
\begin{align*}
\mathrm{(vii)}\hspace{0.5in}&\text{Let}\; S\equiv\pm\{3, 4, 5, 7, 9, 11, 13, 15, 16, 17, 22, 24\} \pmod {50},\; \text{and} \\
&T\equiv\pm\{1, 4, 6, 9, 11, 13, 15, 17, 19, 21, 22, 24\} \pmod{50}.\; \text{Then}\\
&p(S,n)=p(T,n-3),\; \text{for all}\; n\geq 3.
\end{align*}
\begin{align*}
\mathrm{(viii)}\hspace{0.5in}&\text{Let}\; S\equiv\pm\{3, 4, 5, 7, 8, 9, 13, 17, 18, 19, 22, 23\} \pmod {50},\; \text{and} \\
&T\equiv\pm\{1, 4, 5, 9, 12, 13, 15, 18, 19, 21, 22, 23\} \pmod{50}.\; \text{Then}\\
&p(S,n)=p(T,n-3),\; \text{for all}\; n\geq 3.
\end{align*}
\begin{align*}
\mathrm{(ix)}\hspace{0.5in}&\text{Let}\; S\equiv\pm\{2, 5, 7, 8, 9, 11, 12, 14, 15, 17, 19, 21\} \pmod {50},\; \text{and} \\
&T\equiv\pm\{2, 3, 5, 7, 12, 13, 14, 17, 19, 21, 22, 23\} \pmod{50}.\; \text{Then}\\
&p(S,n)=p(T,n-2),\; \text{for all}\; n\geq 2.
\end{align*}
\begin{align*}
\mathrm{(x)}\hspace{0.5in}&\text{Let}\; S\equiv\pm\{3, 4, 5, 6, 7, 9, 15, 17, 18, 21, 23, 24\} \pmod {50},\; \text{and} \\
&T\equiv\pm\{1, 3, 6, 9, 11, 14, 15, 18, 19, 21, 23, 24\} \pmod{50}.\; \text{Then}\\
&p(S,n)=p(T,n-3),\; \text{for all}\; n\geq 3.
\end{align*}
\end{theorem}
\begin{proof}
We prove (i) by starting to write it in its equivalent form
\begin{align}
&[2,9, 15, 21:50]-q[3, 7, 13, 18:50]\notag\\&= [1, 2, 3, 5, 7, 9,
11, 12, 13, 15, 16, 17, 18, 21, 22, 23:50].\label{b1}
\end{align}
As in the previous proof, we use the first identity in Theorem
\ref{50.1} which we express in the form
\begin{align*}
&[2,12,21,23:50]-q[3,8,19,22:50]\\&
=[1,2,3,5,7,8,12,13,15,17,18,19,20,21,22,23:50].
\end{align*}
Upon comparing these last two equations, we see that it suffices to
show
\begin{align*}
&[8,19,20:50]\big\{[2,9, 15, 21:50]-q[3, 7, 13,
18:50]\big\}\\\notag
=&[9,11,16:50]\big\{[2,12,21,23:50]-q[3,8,19,22:50]\big\},
\end{align*}
which after rearrangement takes the form
\begin{align}
&[2,9,21:50]\big\{[11,12,16,23:50]-[8,15,19,20:50]\big\}
\notag\\&=q[3,8,19:50]\big\{[9,11,16,22:50]-[7,13,18,20:50]\big\}.\label{b2}
\end{align}
We employ \eqref{four} twice with the sets of parameters
$[1,4,12,16,20]$ and $[1,3,10,14,19]$ and with $q$ replaced by
$q^{50}$ in both instances, we find that
\begin{equation*}
[11,12,16,23:50]-[8,15,19,20:50]=q^8[3,4,8,19:50]\; \text{and}
\end{equation*}
\begin{equation*}
[9,11,16,22:50]-[7,13,18,20:50]=q^7[2,4,9,21:50].
\end{equation*}
These two equations clearly imply \eqref{b2} and so the proof of
(i) is complete.
\end{proof}

\section{Modulus $M = 52$}
There is only one  class  of identities under multiplication
by the group $U(52)$. \\
\begin{theorem}\label{52.1}

\begin{align*}
\mathrm{(i)}\hspace{0.5in}&\text{Let}\; S\equiv\pm\{1, 3, 4, 5, 7, 16, 18, 19, 21, 22, 23, 25\} \pmod {52},\; \text{and} \\
&T\equiv\pm\{1, 2, 3, 7, 11, 15, 18, 19, 20, 23, 24, 25\} \pmod{52}.\; \text{Then}\\
&p(S,n)=p(T,n-1),\; \text{for all}\; n\geq 1.
\end{align*}
\begin{align*}
\mathrm{(ii)}\hspace{0.5in}&\text{Let}\; S\equiv\pm\{2, 3, 4, 5, 9, 11, 12, 14, 15, 17, 21, 23\} \pmod {52},\; \text{and} \\
&T\equiv\pm\{2, 3, 5, 6, 7, 8, 9, 17, 19, 20, 21, 23\}\pmod{52}.\; \text{Then}\\
&p(S,n)=p(T,n),\; \text{for all}\; n\neq 4.
\end{align*}
\begin{align*}
\mathrm{(iii)}\hspace{0.5in}&\text{Let}\; S\equiv\pm\{3, 4, 5, 9, 10, 11, 14, 15, 16, 17, 21, 23\} \pmod {52},\; \text{and} \\
&T\equiv\pm\{1, 5, 6, 9, 11, 14, 15, 17, 20, 21, 24, 25\}\pmod{52}.\; \text{Then}\\
&p(S,n)=p(T,n-3),\; \text{for all}\; n\geq 3.
\end{align*}
\begin{align*}
\mathrm{(iv)}\hspace{0.5in}&\text{Let}\; S\equiv\pm\{2, 3, 5, 7, 8, 9, 17, 19, 21, 22, 23, 24\} \pmod {52},\; \text{and} \\
&T\equiv\pm\{1, 3, 5, 7, 12, 14, 16, 19, 21, 22, 23, 25\}\pmod{52}.\; \text{Then}\\
&p(S,n)=p(T,n-2),\; \text{for all}\; n\geq 2.
\end{align*}
\begin{align*}
\mathrm{(v)}\hspace{0.5in}&\text{Let}\; S\equiv\pm\{1, 6, 7, 9, 10, 11, 12, 15, 16, 17, 19, 25\} \pmod {52},\; \text{and} \\
&T\equiv\pm\{1, 5, 6, 8, 9, 11, 15, 17, 18, 21, 24, 25\}\pmod{52}.\; \text{Then}\\
&p(S,n)=p(T,n-1),\; \text{for all}\; n\geq 1.
\end{align*}
\begin{align*}
\mathrm{(vi)}\hspace{0.5in}&\text{Let}\; S\equiv\pm\{1, 4, 7, 9, 10, 11, 12, 15, 17, 19, 22, 25\} \pmod {52},\; \text{and} \\
&T\equiv\pm\{1, 3, 7, 8, 10, 11, 15, 18, 19, 20, 23, 25\}\pmod{52}.\; \text{Then}\\
&p(S,n)=p(T,n-1),\; \text{for all}\; n\geq 1.
\end{align*}
\end{theorem}
\begin{proof}
The identities (i)---(vi) follow from \eqref{four2} with the
following choice of parameters
\begin{equation*}
[1, 2, 4, 12, 23],\;[1, 3, 7, 10, 12],\;[1, 5, 2, 7, 11],\; [1, 4,
2, 6, 9],\;[1, 3, 2, 8, 11],\; [1, 2, 5, 19, 20]
\end{equation*}
and with $q$ replaced by $q^{26}$ in each instance.
\end{proof}
\section{Modulus $M = 54$}
There are four  classes  of identities under multiplication
by the group $U(54)$. \\
\begin{theorem}\label{54.1}

\begin{align*}
\mathrm{(i)}\hspace{0.5in}&\text{Let}\; S\equiv\pm\{1, 4, 5, 7, 11, 13, 16, 17, 19, 20, 22, 23\} \pmod {54},\; \text{and} \\
&T\equiv\pm\{1, 3, 4, 10, 11, 15, 16, 17, 19, 21, 23, 26\} \pmod{54}.\; \text{Then}\\
&p(S,n)=p(T,n-1),\; \text{for all}\; n\geq 1.
\end{align*}
\begin{align*}
\mathrm{(ii)}\hspace{0.5in}&\text{Let}\; S\equiv\pm\{1, 2, 5, 7, 8, 11, 13, 19, 20, 23, 25, 26\} \pmod {54},\; \text{and} \\
&T\equiv\pm\{1, 3, 4, 5, 7, 13, 15, 20, 21, 22, 23, 26\}
\pmod{54}.\; \text{Then}\\
&p(S,n)=p(T,n),\; \text{for all}\; n\neq 2.
\end{align*}
\begin{align*}
\mathrm{(iii)}\hspace{0.5in}&\text{Let}\; S\equiv\pm\{1, 4, 5, 7,
8,11,17, 19, 22, 23, 25, 26\}
\pmod {54},\; \text{and} \\
&T\equiv\pm\{1, 3, 4, 7, 11, 15, 16, 20, 21, 23, 25, 26\}
\pmod{54}.\; \text{Then}\\
&p(S,n)=p(T,n-1),\; \text{for all}\; n\geq 1.
\end{align*}
\begin{align*}
\mathrm{(iv)}\hspace{0.5in}&\text{Let}\; S\equiv\pm\{2, 3, 7, 10, 11, 13,
14, 15, 16, 17, 21, 25\}
\pmod {54},\; \text{and} \\
&T\equiv\pm\{1, 4, 7, 10, 11, 13, 14, 17, 19, 23, 25, 26\}
\pmod{54}.\; \text{Then}\\
&p(S,n)=p(T,n-2),\; \text{for all}\; n\geq 2.
\end{align*}
\begin{align*}
\mathrm{(v)}\hspace{0.5in}&\text{Let}\; S\equiv\pm\{2, 5, 7, 8, 10, 11,
13, 16, 17, 19, 23, 25\}
\pmod {54},\; \text{and} \\
&T\equiv\pm\{2, 3, 5, 8, 13, 14, 15, 19, 21, 22, 23, 25\}
\pmod{54}.\; \text{Then}\\
&p(S,n)=p(T,n-2),\; \text{for all}\; n\geq 2.
\end{align*}
\begin{align*}
\mathrm{(vi)}\hspace{0.5in}&\text{Let}\; S\equiv\pm\{1, 2, 4, 5, 11, 13,
14, 16, 17, 19, 23, 25\}
\pmod {54},\; \text{and} \\
&T\equiv\pm\{1, 2, 3, 8, 10, 13, 14, 15, 17, 19, 21, 25\}
\pmod{54}.\; \text{Then}\\
&p(S,n)=p(T,n),\; \text{for all}\; n\neq 3.
\end{align*}
\begin{align*}
\mathrm{(vii)}\hspace{0.5in}&\text{Let}\; S\equiv\pm\{1, 3, 5, 7, 8, 15,
17, 19, 20, 21, 22, 26\}
\pmod {54},\; \text{and} \\
&T\equiv\pm\{1, 2, 5, 7, 13, 14, 17, 19, 20, 22, 23, 25\}
\pmod{54}.\; \text{Then}\\
&p(S,n)=p(T,n-1),\; \text{for all}\; n\geq 1.
\end{align*}
\begin{align*}
\mathrm{(viii)}\hspace{0.5in}&\text{Let}\; S\equiv\pm\{3, 4, 5, 10, 11,
13, 14, 15, 16, 17, 21, 23\}
\pmod {54},\; \text{and} \\
&T\equiv\pm\{1, 5, 7, 10, 11, 13, 16, 17, 20, 23, 25, 26\}
\pmod{54}.\; \text{Then}\\
&p(S,n)=p(T,n-3),\; \text{for all}\; n\geq 3.
\end{align*}
\begin{align*}
\mathrm{(ix)}\hspace{0.5in}&\text{Let}\; S\equiv\pm\{1, 5, 7, 8, 10, 11,
13, 14, 17, 19, 22, 25\}
\pmod {54},\; \text{and} \\
&T\equiv\pm\{2, 3, 5, 7, 8, 11, 15, 19, 20, 21, 22, 25\}
\pmod{54}.\; \text{Then}\\
&p(S,n)=p(T,n),\; \text{for all}\; n\neq 1.
\end{align*}
\end{theorem}
\begin{proof}
The identities (i)---(ix) follow from \eqref{four2} with the
choice of parameters $[1, 2, 5, 16, 22]$,\;\\$[1, 2, 4, 8, 9]$,\;
$[1, 3, 2, 6, 9],\;[1, 3, 4, 18, 20],\;[1, 5, 3, 8, 11],\;[1, 2,
5, 7, 13],\; [1, 2, 4, 14, 23]$,\\ $[1, 5, 2, 7, 12],\; [1, 3, 4,
9, 11]$ and with $q$ replaced by $q^{27}$ in each instance.
\end{proof}
\begin{theorem}\label{54.2}

\begin{align*}
\mathrm{(i)}\hspace{0.5in}&\text{Let}\; S\equiv\pm\{1, 5, 6, 9, 10, 11, 13, 14, 21, 22, 23, 25\} \pmod {54},\; \text{and} \\
&T\equiv\pm\{1, 4, 5, 9, 10, 11, 17, 19, 21, 22, 23, 24\} \pmod{54}.\; \text{Then}\\
&p(S,n)=p(T,n-1),\; \text{for all}\; n\geq 1.
\end{align*}
\begin{align*}
\mathrm{(ii)}\hspace{0.5in}&\text{Let}\; S\equiv\pm\{1, 2, 3, 4, 5, 7, 9, 12, 13, 20, 23, 25\} \pmod {54},\; \text{and} \\
&T\equiv\pm\{1, 2, 3, 4, 5, 7, 9, 11, 16, 17, 24, 25\} \pmod{54}.\; \text{Then}\\
&p(S,n)=p(T,n),\; \text{for all}\; n\neq 11.
\end{align*}
\begin{align*}
\mathrm{(iii)}\hspace{0.5in}&\text{Let}\; S\equiv\pm\{1, 7, 8, 9, 10, 12, 13, 15, 16, 17, 19, 23\} \pmod {54},\; \text{and} \\
&T\equiv\pm\{1, 6, 7, 8, 9, 11, 15, 16, 19, 23, 25, 26\} \pmod{54}.\; \text{Then}\\
&p(S,n)=p(T,n-1),\; \text{for all}\; n\geq 1.
\end{align*}
\begin{align*}
\mathrm{(iv)}\hspace{0.5in}&\text{Let}\; S\equiv\pm\{1, 2, 6, 7, 9, 10, 11, 13, 15, 17, 25, 26\} \pmod {54},\; \text{and} \\
&T\equiv\pm\{1, 2, 5, 8, 9, 11, 12, 13, 15, 17, 19, 26\} \pmod{54}.\; \text{Then}\\
&p(S,n)=p(T,n),\; \text{for all}\; n\neq 5.
\end{align*}
\begin{align*}
\mathrm{(v)}\hspace{0.5in}&\text{Let}\; S\equiv\pm\{2, 3, 5, 9, 11, 12, 13, 16, 19, 22, 23, 25\} \pmod {54},\; \text{and} \\
&T\equiv\pm\{1, 3, 7, 9, 11, 13, 16, 19, 20, 22, 24, 25\} \pmod{54}.\; \text{Then}\\
&p(S,n)=p(T,n-2),\; \text{for all}\; n\geq 2.
\end{align*}
\begin{align*}
\mathrm{(vi)}\hspace{0.5in}&\text{Let}\; S\equiv\pm\{4, 5, 6, 7, 8, 9, 13, 17, 21, 22, 23, 25\} \pmod {54},\; \text{and} \\
&T\equiv\pm\{1, 4, 8, 9, 13, 14, 17, 19, 21, 23, 24, 25\} \pmod{54}.\; \text{Then}\\
&p(S,n)=p(T,n-4),\; \text{for all}\; n\geq 4.
\end{align*}
\begin{align*}
\mathrm{(vii)}\hspace{0.5in}&\text{Let}\; S\equiv\pm\{4, 5, 6, 7, 9, 11, 13, 14, 19, 21, 23, 26\} \pmod {54},\; \text{and} \\
&T\equiv\pm\{1, 5, 7, 9, 13, 14, 17, 19, 21, 22, 24, 26\} \pmod{54}.\; \text{Then}\\
&p(S,n)=p(T,n-4),\; \text{for all}\; n\geq 4.
\end{align*}
\begin{align*}
\mathrm{(viii)}\hspace{0.5in}&\text{Let}\; S\equiv\pm\{3, 5, 7, 9, 11, 12, 13, 14, 16, 17, 20, 23\} \pmod {54},\; \text{and} \\
&T\equiv\pm\{2, 3, 7, 9, 11, 14, 17, 19, 20, 23, 24, 25\} \pmod{54}.\; \text{Then}\\
&p(S,n)=p(T,n-3),\; \text{for all}\; n\geq 3.
\end{align*}
\begin{align*}
\mathrm{(ix)}\hspace{0.5in}&\text{Let}\; S\equiv\pm\{5, 6, 7, 8, 9, 10, 11, 15, 17, 19, 20, 25\} \pmod {54},\; \text{and} \\
&T\equiv\pm\{1, 5, 9, 10, 12, 15, 17, 19, 20, 23, 25, 26\} \pmod{54}.\; \text{Then}\\
&p(S,n)=p(T,n-5),\; \text{for all}\; n\geq 5.
\end{align*}
\end{theorem}
\begin{proof}
We prove (ii). We apply \eqref{qp} with $q$ replaced by $q^9$ and
$x$ is replaced by $q^2$, we find that
\begin{align}
&[3,24,24:54]-q^2[15,12,6:54]\notag\\
&=[ 2, 3, 5, 7, 9, 11, 12, 13, 15,16,18, 20, 23, 24,
25:54].\label{memo}
\end{align}
We reformulate the identity in (ii) in its equivalent form
\begin{align}
&[12,13,20,23:54]-[11,16,17,24:54]\notag\\&=q^{11}[1, 2, 3, 4, 5,
7, 9, 11, 12, 13, 16, 17, 20, 23, 24, 25:54].\label{memo2}
\end{align}
By \eqref{memo}, the equation \eqref{memo2} is equivalent to
\begin{align}
&q^{11}[1,4,17]\big\{[3,24,24:54]-q^2[6,12,15:54]\big\}\notag\\
&=[15,18]\big\{[12,13,20,23:54]-[11,16,17,24:54]\big\}.\label{memo3}
\end{align}
After regrouping  the the terms in \eqref{memo3}, we  are led to
prove
\begin{align}
&[17,24:54]\big\{[11,15,16,18:54]+q^{11}[1,3,4,24:54]\big\}\notag\\
&=[15,12:54]\big\{[18,13,20,23:54]+q^{13}[1,4,6,17:54]\big\}.\label{memo4}
\end{align}
Next, we apply \eqref{four} twice with the sets of parameters
$[1,2,13,16,17]$, $[1,2,15,19,21]$ and with $q$ replaced by $q^{54}$
in each instance to find that
\begin{align*}
&[11,15,16,18:54]+q^{11}[1,3,4,24:54]=[12,14,15,19:54] \;\text{and}\\
&[18,13,20,23]+q^{13}[1,4,6,17:54]=[14,17,19,24].
\end{align*}
These two identities clearly imply \eqref{memo4}. Hence, the proof
of (ii) is complete.
\end{proof}
\begin{theorem}\label{54.3}

\begin{align*}
\mathrm{(i)}\hspace{0.5in}&\text{Let}\; S\equiv\pm\{1, 3, 5, 7, 11, 13, 15, 16, 20, 24, 25, 26\} \pmod {54},\; \text{and} \\
&T\equiv\pm\{1, 2, 5, 9, 12, 13, 15, 19, 20, 23, 25, 26\} \pmod{54}.\; \text{Then}\\
&p(S,n)=p(T,n-1),\; \text{for all}\; n\geq 1.
\end{align*}
\begin{align*}
\mathrm{(ii)}\hspace{0.5in}&\text{Let}\; S\equiv\pm\{5, 6, 7, 8, 9, 10, 11, 13, 17, 21, 22, 25\} \pmod {54},\; \text{and} \\
&T\equiv\pm\{1, 5, 8, 11, 12, 15, 17, 19, 21, 22, 25, 26\} \pmod{54}.\; \text{Then}\\
&p(S,n)=p(T,n-5),\; \text{for all}\; n\geq 5.
\end{align*}
\begin{align*}
\mathrm{(iii)}\hspace{0.5in}&\text{Let}\; S\equiv\pm\{3, 4, 5, 6, 7, 13, 17, 19, 20, 21, 22, 23\} \pmod {54},\; \text{and} \\
&T\equiv\pm\{1, 3, 7, 9, 13, 14, 17, 19, 20, 22, 24, 25\} \pmod{54}.\; \text{Then}\\
&p(S,n)=p(T,n-3),\; \text{for all}\; n\geq 3.
\end{align*}
\begin{align*}
\mathrm{(iv)}\hspace{0.5in}&\text{Let}\; S\equiv\pm\{1, 3, 4, 5, 7, 9, 11, 16, 17, 19, 22, 24\} \pmod {54},\; \text{and} \\
&T\equiv\pm\{1, 3, 4, 5, 6, 11, 13, 14, 16, 19, 21, 23\} \pmod{54}.\; \text{Then}\\
&p(S,n)=p(T,n),\; \text{for all}\; n\neq 6.
\end{align*}
\begin{align*}
\mathrm{(v)}\hspace{0.5in}&\text{Let}\; S\equiv\pm\{1, 7, 8, 10, 11, 12, 13, 14, 15, 17, 19, 21\} \pmod {54},\; \text{and} \\
&T\equiv\pm\{1, 6, 7, 9, 10, 11, 13, 14, 21, 23, 25, 26\} \pmod{54}.\; \text{Then}\\
&p(S,n)=p(T,n-1),\; \text{for all}\; n\geq 1.
\end{align*}
\begin{align*}
\mathrm{(vi)}\hspace{0.5in}&\text{Let}\; S\equiv\pm\{1, 5, 7, 9, 10, 12, 13, 15, 16, 17, 20, 23\} \pmod {54},\; \text{and} \\
&T\equiv\pm\{2, 3, 5, 7, 10, 11, 15, 16, 17, 23, 24, 25\} \pmod{54}.\; \text{Then}\\
&p(S,n)=p(T,n),\; \text{for all}\; n\neq 1.
\end{align*}
\begin{align*}
\mathrm{(vii)}\hspace{0.5in}&\text{Let}\; S\equiv\pm\{2, 5, 8, 9, 11, 12, 13, 15, 16, 17, 19, 23\} \pmod {54},\; \text{and} \\
&T\equiv\pm\{2, 3, 7, 8, 11, 13, 15, 19, 20, 23, 24, 25\} \pmod{54}.\; \text{Then}\\
&p(S,n)=p(T,n-2),\; \text{for all}\; n\geq 2.
\end{align*}
\begin{align*}
\mathrm{(viii)}\hspace{0.5in}&\text{Let}\; S\equiv\pm\{4, 5, 6, 7, 8, 9, 11, 19, 21, 23, 25, 26\} \pmod {54},\; \text{and} \\
&T\equiv\pm\{1, 4, 7, 10, 12, 15, 17, 19, 21, 23, 25, 26\} \pmod{54}.\; \text{Then}\\
&p(S,n)=p(T,n-4),\; \text{for all}\; n\geq 4.
\end{align*}
\begin{align*}
\mathrm{(ix)}\hspace{0.5in}&\text{Let}\; S\equiv\pm\{1, 2, 3, 4, 9, 11, 14, 17, 19, 23, 24, 25\} \pmod {54},\; \text{and} \\
&T\equiv\pm\{1, 2, 3, 5, 6, 13, 14, 17, 21, 22, 23, 25\} \pmod{54}.\; \text{Then}\\
&p(S,n)=p(T,n),\; \text{for all}\; n\neq 4.
\end{align*}
\end{theorem}
\begin{proof}
We prove (i) by using the same identity \eqref{memo} that we
employed in our previous proof namely,
\begin{align}
&[3,24,24:54]-q^2[15,12,6:54]\notag\\
&=[ 2, 3, 5, 7, 9, 11, 12, 13, 15,16,18, 20, 23, 24,
25:54].\label{memo6}
\end{align}
We start by writing (i) in its equivalent form
\begin{align}
&[2,9,12,19,23:54]-q[3,7,11,16,24:54]\notag\\
&=[1, 2, 3, 5, 7, 9,
11, 12, 13, 15, 16, 19, 20, 23, 24, 25, 26:54].\label{memo7}
\end{align}
By \eqref{memo6}, we see that \eqref{memo7} is equivalent to
\begin{align}
&[18:54]\big\{[2,9,12,19,23:54]-q[3,7,11,16,24:54]\big\}\notag\\
&=[1,19,26:54]\big\{[3,24,24:54]-q^2[6,12,15:54]\big\}.\label{memo8}
\end{align}
After rearrangement of the terms in above sum, we are led to prove
\begin{align}
&[3,24:54]\big\{[1,19,24,26:54]+q[7,11,16,18:54]\big\}\notag\\
&=[12,19:54]\big\{[2,9,18,23:54]+q^2[1,6,15,26:54]\big\}.\label{memo9}
\end{align}
The equation \eqref{memo9} follows from the pair of identities
\begin{align*}
&[1,19,24,26:54]+q[7,11,16,18:54]=[8,12,17,19:54] \;\text{and}\\
&[2,9,18,23:54]+q^2[1,6,15,26:54]=[3,8,17,24:54],
\end{align*}
which are obtained by \eqref{four} by the choice of parameters
$[1,8,9,20,25],\;[1,2,4,10,19]$ and $q$ replaced by $q^{54}$ in
each instance.
\end{proof}
\begin{theorem}\label{54.4}

\begin{align*}
\mathrm{(i)}\hspace{0.5in}&\text{Let}\; S\equiv\pm\{1, 3, 5, 7, 11, 12, 13, 20, 22, 23, 24, 25\} \pmod {54},\; \text{and} \\
&T\equiv\pm\{1, 2, 5, 9, 11, 12, 17, 19, 22, 23, 24, 25\} \pmod{54}.\; \text{Then}\\
&p(S,n)=p(T,n-1),\; \text{for all}\; n\geq 1.
\end{align*}
\begin{align*}
\mathrm{(ii)}\hspace{0.5in}&\text{Let}\; S\equiv\pm\{1, 2, 5, 6, 7, 9, 10, 12, 13, 17, 23, 25\} \pmod {54},\; \text{and} \\
&T\equiv\pm\{1, 2, 5, 6, 7, 8, 11, 12, 15, 17, 19, 25\} \pmod{54}.\; \text{Then}\\
&p(S,n)=p(T,n),\; \text{for all}\; n\neq 8.
\end{align*}
\begin{align*}
\mathrm{(iii)}\hspace{0.5in}&\text{Let}\; S\equiv\pm\{1, 6, 7, 8, 9, 11, 13, 14, 19, 23, 24, 25\} \pmod {54},\; \text{and} \\
&T\equiv\pm\{1, 5, 6, 7, 8, 13, 17, 19, 21, 22, 23, 24\} \pmod{54}.\; \text{Then}\\
&p(S,n)=p(T,n-1),\; \text{for all}\; n\geq 1.
\end{align*}
\begin{align*}
\mathrm{(iv)}\hspace{0.5in}&\text{Let}\; S\equiv\pm\{1, 5, 6, 7, 9, 11, 13, 17, 22, 24, 25, 26\} \pmod {54},\; \text{and} \\
&T\equiv\pm\{1, 4, 5, 6, 11, 13, 17, 19, 21, 23, 24, 26\} \pmod{54}.\; \text{Then}\\
&p(S,n)=p(T,n-1),\; \text{for all}\; n\geq 1.
\end{align*}
\begin{align*}
\mathrm{(v)}\hspace{0.5in}&\text{Let}\; S\equiv\pm\{1, 6, 7, 10, 11, 12, 13, 15, 16, 17, 19, 25\} \pmod {54},\; \text{and} \\
&T\equiv\pm\{1, 5, 6, 9, 11, 12, 13, 16, 19, 23, 25, 26\} \pmod{54}.\; \text{Then}\\
&p(S,n)=p(T,n-1),\; \text{for all}\; n\geq 1.
\end{align*}
\begin{align*}
\mathrm{(vi)}\hspace{0.5in}&\text{Let}\; S\equiv\pm\{3, 4, 5, 7, 11, 12, 13, 16, 17, 23, 24, 25\} \pmod {54},\; \text{and} \\
&T\equiv\pm\{1, 4, 7, 9, 12, 13, 17, 19, 20, 23, 24, 25\} \pmod{54}.\; \text{Then}\\
&p(S,n)=p(T,n-3),\; \text{for all}\; n\geq 3.
\end{align*}
\begin{align*}
\mathrm{(vii)}\hspace{0.5in}&\text{Let}\; S\equiv\pm\{1, 5, 7, 9, 11, 12, 13, 14, 16, 17, 19, 24\}\pmod {54},\; \text{and} \\
&T\equiv\pm\{2, 3, 5, 7, 11, 12, 13, 14, 19, 23, 24, 25\} \pmod{54}.\; \text{Then}\\
&p(S,n)=p(T,n),\; \text{for all}\; n\neq 1.
\end{align*}
\begin{align*}
\mathrm{(viii)}\hspace{0.5in}&\text{Let}\; S\equiv\pm\{5, 6, 7, 8, 9, 11, 12, 13, 17, 19, 20, 23\}\pmod {54},\; \text{and} \\
&T\equiv\pm\{1, 6, 7, 11, 12, 15, 17, 19, 20, 23, 25, 26\} \pmod{54}.\; \text{Then}\\
&p(S,n)=p(T,n-5),\; \text{for all}\; n\geq 5.
\end{align*}
\begin{align*}
\mathrm{(ix)}\hspace{0.5in}&\text{Let}\; S\equiv\pm\{4, 5, 6, 7, 9, 10, 11, 17, 19, 23, 24, 25\}\pmod {54},\; \text{and} \\
&T\equiv\pm\{1, 5, 6, 10, 13, 14, 17, 19, 21, 23, 24, 25\} \pmod{54}.\; \text{Then}\\
&p(S,n)=p(T,n-4),\; \text{for all}\; n\geq 4.
\end{align*}
\end{theorem}
\begin{proof}

We prove (i) which is equivalent to
\begin{align}
&[2,9,17,19:54]-q[3,7,13,20:54]\notag\\
&=[1, 2, 3, 5, 7, 9, 11,12, 13, 17, 19, 20, 22, 23, 24,
25:54].\label{mu1}
\end{align}
By \eqref{memo2}
\begin{align}
&[12,13,20,23:54]-[11,16,17,24:54]\notag\\&=q^{11}[1, 2, 3, 4, 5,
7, 9, 11, 12, 13, 16, 17, 20, 23, 24, 25:54].\label{memo22}
\end{align}
By \eqref{memo22}, the equation \eqref{mu1} is equivalent to
\begin{align}
&q^{11}[4,16:54]\big\{[2,9,17,19:54]-q[3,7,13,20:54]\big\}\notag\\
&=[19,22:54]\big\{[12,13,20,23:54]-[11,16,17,24:54]\big\}.\label{mu2}
\end{align}
After regrouping  the terms in \eqref{mu2}, we  are led to
prove
\begin{align}
&[13,20:54]\big\{[12,19,22,23:54]+q^{12}[3,4,7,16:54]\big\}\notag\\
&=[19,16:54]\big\{[11,17,22,24]+q^{11}[2,4,9,17:54]\big\}.\label{mu3}
\end{align}
Next, we apply \eqref{four} twice with the sets of parameters
$[1,4,16,20,23]$, $[1,3,14,18,23]$ and with $q$ replaced by $q^{54}$
in each instance to find that
\begin{align*}
&[12,19,22,23:54]+q^{12}[3,4,7,16:54]=[15,16,19,26:54] \;\text{and}\\
&[11,17,22,24]+q^{11}[2,4,9,17:54]=[13,15,20,26:54].
\end{align*}
These two identities readily imply \eqref{mu3}. Hence, the proof
of (i) is complete.
\end{proof}
\section{Modulus $M = 56$}
There is only one  class  of identities under multiplication
by the group $U(56)$. \\
\begin{theorem}\label{56.1}

\begin{align*}
\mathrm{(i)}\hspace{0.5in}&\text{Let}\; S\equiv\pm\{1, 3, 4, 5, 7, 18, 21, 22, 23, 24, 25, 27\} \pmod {56},\; \text{and} \\
&T\equiv\pm\{1, 2, 3, 7, 11, 17, 20, 21, 24, 25, 26, 27\} \pmod{56}.\; \text{Then}\\
&p(S,n)=p(T,n-1),\; \text{for all}\; n\geq 1.
\end{align*}
\begin{align*}
\mathrm{(ii)}\hspace{0.5in}&\text{Let}\; S\equiv\pm\{2, 3, 7, 9, 10, 12, 13, 15, 16, 19, 21, 25\} \pmod {56},\; \text{and} \\
&T\equiv\pm\{3, 4, 5, 6, 7, 9, 16, 19, 21, 22, 23, 25\}\pmod{56}.\; \text{Then}\\
&p(S,n)=p(T,n),\; \text{for all}\; n\neq 2.
\end{align*}
\begin{align*}
\mathrm{(iii)}\hspace{0.5in}&\text{Let}\; S\equiv\pm\{2, 3, 5, 7, 8, 13, 15, 20, 21, 22, 23, 25\} \pmod {56},\; \text{and} \\
&T\equiv\pm\{1, 5, 7, 8, 10, 12, 13, 15, 18, 21, 23, 27\}\pmod{56}.\; \text{Then}\\
&p(S,n)=p(T,n),\; \text{for all}\; n\neq 1.
\end{align*}
\begin{align*}
\mathrm{(iv)}\hspace{0.5in}&\text{Let}\; S\equiv\pm\{1, 7, 8, 9, 10, 12, 13, 15, 18, 19, 21, 27\} \pmod {56},\; \text{and} \\
&T\equiv\pm\{1, 6, 7, 8, 9, 11, 17, 19, 20, 21, 26, 27\}\pmod{56}.\; \text{Then}\\
&p(S,n)=p(T,n-1),\; \text{for all}\; n\geq 1.
\end{align*}
\begin{align*}
\mathrm{(v)}\hspace{0.5in}&\text{Let}\; S\equiv\pm\{4, 5, 6, 7, 9, 11, 16, 17, 19, 21, 22, 23\} \pmod {56},\; \text{and} \\
&T\equiv\pm\{1, 5, 7, 11, 12, 16, 17, 18, 21, 23, 26, 27\}\pmod{56}.\; \text{Then}\\
&p(S,n)=p(T,n-4),\; \text{for all}\; n\geq 4.
\end{align*}
\begin{align*}
\mathrm{(vi)}\hspace{0.5in}&\text{Let}\; S\equiv\pm\{4, 6, 7, 9, 10, 11, 13, 15, 17, 19, 21, 24\} \pmod {56},\; \text{and} \\
&T\equiv\pm\{2, 3, 7, 11, 13, 15, 17, 20, 21, 24, 25, 26\}\pmod{56}.\; \text{Then}\\
&p(S,n)=p(T,n-4),\; \text{for all}\; n\geq 4.
\end{align*}
\end{theorem}
\begin{proof}
The identities (i)---(vi) follow from \eqref{four2} with the
following choice of parameters
\begin{equation*}
[1, 2, 4, 12, 25],\;[1, 4, 6, 10, 13],\; [1, 3, 4, 9, 11], \;[1, 3,
2, 9, 11],\; [1, 5, 6, 17, 24],\; [1, 7, 3, 10, 11]
\end{equation*}
and with $q$ replaced by $q^{28}$ in each instance.
\end{proof}
\section{Modulus $M = 60$}
There are eight  classes  of identities under multiplication
by the group $U(60)$. The first two identities in Theorem \ref{60.3} below 
were  first given  by Kalvade \cite{KA}. \\
\begin{theorem}\label{60.1}

\begin{align*}
\mathrm{(i)}\hspace{0.5in}&\text{Let}\; S\equiv\pm\{1, 3, 4, 5, 7, 22, 23, 24, 25, 26, 27, 29\}\pmod {60},\; \text{and} \\
&T\equiv\pm\{1, 2, 3, 7, 11, 19, 23, 24, 26, 27, 28, 29\} \pmod{60}.\; \text{Then}\\
&p(S,n)=p(T,n-1),\; \text{for all}\; n\geq 1.
\end{align*}
\begin{align*}
\mathrm{(ii)}\hspace{0.5in}&\text{Let}\; S\equiv\pm\{2, 7, 9, 11, 12, 13, 14, 16, 17, 19, 21, 23\}\pmod {60},\; \text{and} \\
&T\equiv\pm\{2, 5, 7, 9, 11, 12, 19, 21, 23, 25, 26, 28\} \pmod{60}.\; \text{Then}\\
&p(S,n)=p(T,n-2),\; \text{for all}\; n\geq 2.
\end{align*}
\begin{align*}
\mathrm{(iii)}\hspace{0.5in}&\text{Let}\; S\equiv\pm\{2, 3, 5, 11, 13, 14, 16, 17, 19, 24, 25, 27\}\pmod {60},\; \text{and} \\
&T\equiv\pm\{1, 3, 8, 11, 13, 14, 17, 19, 22, 24, 27, 29\} \pmod{60}.\; \text{Then}\\
&p(S,n)=p(T,n-2),\; \text{for all}\; n\geq 2.
\end{align*}
\begin{align*}
\mathrm{(iv)}\hspace{0.5in}&\text{Let}\; S\equiv\pm\{1, 5, 8, 9, 12, 13, 14, 17, 21, 22, 25, 29\}\pmod {60},\; \text{and} \\
&T\equiv\pm\{1, 4, 7, 9, 12, 13, 17, 21, 22, 23, 26, 29\} \pmod{60}.\; \text{Then}\\
&p(S,n)=p(T,n-1),\; \text{for all}\; n\geq 1.
\end{align*}
\end{theorem}
\begin{proof}
The identities (i)--(iv) follow from \eqref{four2} with the
following sets of parameters
\begin{equation*}
[1,2,4,12,27],\;[1,5,3,10,14],\;[1,3,4,17,23],\;[1,2,6,22,23],
\end{equation*}
and with $q$ replaced by $q^{30}$ in each instance.
\end{proof}
\begin{theorem}\label{60.2}

\begin{align*}
\mathrm{(i)}\hspace{0.5in}&\text{Let}\; S\equiv\pm\{1, 5, 6, 8, 9, 13, 17, 21, 22, 25, 28, 29\}\pmod {60},\; \text{and} \\
&T\equiv\pm\{1, 4, 5, 7, 13, 16, 17, 22, 23, 25, 26, 29\} \pmod{60}.\; \text{Then}\\
&p(S,n)=p(T,n-1),\; \text{for all}\; n\geq 1.
\end{align*}
\begin{align*}
\mathrm{(ii)}\hspace{0.5in}&\text{Let}\; S\equiv\pm\{1, 2, 5, 7, 8, 11, 19, 23, 25, 26, 28, 29\}\pmod {60},\; \text{and} \\
&T\equiv\pm\{1, 3, 4, 5, 7, 16, 18, 23, 25, 26, 27, 29\}\pmod{60}.\; \text{Then}\\
&p(S,n)=p(T,n),\; \text{for all}\; n\neq 2.
\end{align*}
\begin{align*}
\mathrm{(iii)}\hspace{0.5in}&\text{Let}\; S\equiv\pm\{2, 4, 5, 7, 11, 13, 14, 16, 17, 19, 23, 25\}\pmod {60},\; \text{and} \\
&T\equiv\pm\{2, 5, 6, 7, 8, 9, 11, 19, 21, 23, 25, 28\} \pmod{60}.\; \text{Then}\\
&p(S,n)=p(T,n),\; \text{for all}\; n\neq 4.
\end{align*}
\begin{align*}
\mathrm{(iv)}\hspace{0.5in}&\text{Let}\; S\equiv\pm\{3, 4, 5, 11, 13, 14, 16, 17, 18, 19, 25, 27\}\pmod {60},\; \text{and} \\
&T\equiv\pm\{1, 5, 8, 11, 13, 14, 17, 19, 22, 25, 28, 29\} \pmod{60}.\; \text{Then}\\
&p(S,n)=p(T,n-3),\; \text{for all}\; n\geq 3.
\end{align*}
\end{theorem}
\begin{proof}
The identities (i)--(iv) follow from \eqref{four2} with the
following sets of parameters
\begin{equation*}
[1,3,2,7,10],\;[1,2,4,8,9],\;[1,3,7,12,14],\;[1,5,2,7,13],
\end{equation*}
and with $q$ replaced by $q^{30}$ in each instance.
\end{proof}

\begin{theorem}\label{60.3}

\begin{align*}
\mathrm{(i)}\hspace{0.5in}&\text{Let}\; S\equiv\pm\{1, 6, 7, 11, 12, 13, 16, 17, 18, 19, 23, 29\}\pmod {60},\; \text{and} \\
&T\equiv\pm\{1, 5, 6, 11, 12, 13, 17, 18, 19, 25, 28, 29\} \pmod{60}.\; \text{Then}\\
&p(S,n)=p(T,n-1),\; \text{for all}\; n\geq 1.
\end{align*}
\begin{align*}
\mathrm{(ii)}\hspace{0.5in}&\text{Let}\; S\equiv\pm\{1, 6, 7, 8, 11, 13, 17, 18, 19, 23, 24, 29\}\pmod {60},\; \text{and} \\
&T\equiv\pm\{1, 5, 6, 7, 13, 16, 17, 18, 23, 24, 25, 29\}\pmod{60}.\; \text{Then}\\
&p(S,n)=p(T,n-1),\; \text{for all}\; n\geq 1.
\end{align*}
\begin{align*}
\mathrm{(iii)}\hspace{0.5in}&\text{Let}\; S\equiv\pm\{1, 5, 6, 7, 8, 11, 12, 18, 19, 23, 25, 29\}\pmod {60},\; \text{and} \\
&T\equiv\pm\{1, 4, 6, 7, 11, 12, 13, 17, 18, 19, 23, 29\} \pmod{60}.\; \text{Then}\\
&p(S,n)=p(T,n),\; \text{for all}\; n\neq 4.
\end{align*}
\begin{align*}
\mathrm{(iv)}\hspace{0.5in}&\text{Let}\; S\equiv\pm\{4, 5, 6, 7, 11, 13, 17, 18, 19, 23, 24, 25\}\pmod {60},\; \text{and} \\
&T\equiv\pm\{1, 6, 7, 11, 13, 17, 18, 19, 23, 24, 28, 29\} \pmod{60}.\; \text{Then}\\
&p(S,n)=p(T,n-4),\; \text{for all}\; n\geq 4.
\end{align*}
\end{theorem}
\begin{proof}
The identities (i)--(iv) follow from \eqref{four2} with the
following sets of parameters
\begin{equation*}
[1,3,2,8,13],\;[1,2,7,21,24],\;[1,2,6,12,13],\;[1,5,6,18,25],
\end{equation*}
and with $q$ replaced by $q^{30}$ in each instance.
\end{proof}
\begin{theorem}\label{60.4}

\begin{align*}
\mathrm{(i)}\hspace{0.5in}&\text{Let}\; S\equiv\pm\{1, 5, 6, 9, 13, 14, 16, 17, 20, 21, 25, 29\}\pmod {60},\; \text{and} \\
&T\equiv\pm\{1, 4, 5, 11, 13, 14, 17, 19, 20, 25, 26, 29\} \pmod{60}.\; \text{Then}\\
&p(S,n)=p(T,n-1),\; \text{for all}\; n\geq 1.
\end{align*}
\begin{align*}
\mathrm{(ii)}\hspace{0.5in}&\text{Let}\; S\equiv\pm\{1, 3, 5, 7, 8, 18, 20, 22, 23, 25, 27, 29\}\pmod {60},\; \text{and} \\
&T\equiv\pm\{1, 2, 5, 7, 13, 17, 20, 22, 23, 25, 28, 29\}\pmod{60}.\; \text{Then}\\
&p(S,n)=p(T,n-1),\; \text{for all}\; n\geq 1.
\end{align*}
\begin{align*}
\mathrm{(iii)}\hspace{0.5in}&\text{Let}\; S\equiv\pm\{4, 5, 6, 7, 9, 11, 19, 20, 21, 23, 25, 26\}\pmod {60},\; \text{and} \\
&T\equiv\pm\{1, 5, 7, 11, 14, 16, 19, 20, 23, 25, 26, 29\} \pmod{60}.\; \text{Then}\\
&p(S,n)=p(T,n-4),\; \text{for all}\; n\geq 4.
\end{align*}
\begin{align*}
\mathrm{(iv)}\hspace{0.5in}&\text{Let}\; S\equiv\pm\{2, 5, 7, 8, 11, 13, 17, 19, 20, 22, 23, 25\}\pmod {60},\; \text{and} \\
&T\equiv\pm\{2, 3, 5, 11, 13, 17, 18, 19, 20, 25, 27, 28\} \pmod{60}.\; \text{Then}\\
&p(S,n)=p(T,n-2),\; \text{for all}\; n\geq 2.
\end{align*}
\end{theorem}
\begin{proof}
The identities (i)--(iv) follow from \eqref{four2} with the
following sets of parameters
\begin{equation*}
[1,2,6,20,23],\;[1,2,4,14,26],\;[1,5,6,17,26],\;[1,3,6,19,25],
\end{equation*}
and with $q$ replaced by $q^{30}$ in each instance.
\end{proof}
\begin{theorem}\label{60.5}

\begin{align*}
\mathrm{(i)}\hspace{0.5in}&\text{Let}\; S\equiv\pm\{1, 7, 8, 10, 11, 13, 14, 17, 19, 20, 23, 29\}\pmod {60},\; \text{and} \\
&T\equiv\pm\{1, 6, 7, 9, 10, 13, 17, 20, 21, 23, 28, 29\} \pmod{60}.\; \text{Then}\\
&p(S,n)=p(T,n-1),\; \text{for all}\; n\geq 1.
\end{align*}
\begin{align*}
\mathrm{(ii)}\hspace{0.5in}&\text{Let}\; S\equiv\pm\{1, 4, 7, 10, 11, 13, 17, 19, 20, 22, 23, 29\}\pmod {60},\; \text{and} \\
&T\equiv\pm\{1, 3, 7, 10, 11, 16, 18, 19, 20, 23, 27, 29\}\pmod{60}.\; \text{Then}\\
&p(S,n)=p(T,n-1),\; \text{for all}\; n\geq 1.
\end{align*}
\begin{align*}
\mathrm{(iii)}\hspace{0.5in}&\text{Let}\; S\equiv\pm\{6, 7, 8, 9, 10, 11, 13, 17, 19, 20, 21, 23\}\pmod {60},\; \text{and} \\
&T\equiv\pm\{1, 7, 10, 11, 13, 17, 19, 20, 23, 26, 28, 29\} \pmod{60}.\; \text{Then}\\
&p(S,n)=p(T,n-6),\; \text{for all}\; n\geq 6.
\end{align*}
\begin{align*}
\mathrm{(iv)}\hspace{0.5in}&\text{Let}\; S\equiv\pm\{1, 3, 4, 10, 11, 13, 17, 18, 19, 20, 27, 29\}\pmod {60},\; \text{and} \\
&T\equiv\pm\{1, 2, 7, 10, 11, 13, 16, 17, 19, 20, 23, 29\} \pmod{60}.\; \text{Then}\\
&p(S,n)=p(T,n),\; \text{for all}\; n\neq 2.
\end{align*}
\end{theorem}
\begin{proof}
The identities (i)--(iv) follow from \eqref{four2} with the
following sets of parameters
\begin{equation*}
[1,3,2,9,12],\;[1,2,5,19,24],\;[1,8,2,10,12],\;[1,2,4,8,14],
\end{equation*}
and with $q$ replaced by $q^{30}$ in each instance.
\end{proof}
\begin{theorem}\label{60.6}

\begin{align*}
\mathrm{(i)}\hspace{0.5in}&\text{Let}\; S\equiv\pm\{1, 8, 10, 11, 12, 13, 14, 15, 17, 19, 21, 23\}\pmod {60},\; \text{and} \\
&T\equiv\pm\{1, 7, 9, 10, 11, 12, 13, 15, 23, 26, 28, 29\} \pmod{60}.\; \text{Then}\\
&p(S,n)=p(T,n-1),\; \text{for all}\; n\geq 1.
\end{align*}
\begin{align*}
\mathrm{(ii)}\hspace{0.5in}&\text{Let}\; S\equiv\pm\{2, 3, 7, 10, 11, 15, 16, 17, 19, 23, 24, 29\}\pmod {60},\; \text{and} \\
&T\equiv\pm\{1, 4, 7, 10, 13, 15, 17, 19, 22, 24, 27, 29\}\pmod{60}.\; \text{Then}\\
&p(S,n)=p(T,n-2),\; \text{for all}\; n\geq 2.
\end{align*}
\begin{align*}
\mathrm{(iii)}\hspace{0.5in}&\text{Let}\; S\equiv\pm\{1, 3, 4, 10, 11, 13, 15, 17, 22, 23, 24, 29\}\pmod {60},\; \text{and} \\
&T\equiv\pm\{1, 2, 7, 10, 11, 13, 15, 16, 19, 23, 24, 27\} \pmod{60}.\; \text{Then}\\
&p(S,n)=p(T,n),\; \text{for all}\; n\neq 2.
\end{align*}
\begin{align*}
\mathrm{(iv)}\hspace{0.5in}&\text{Let}\; S\equiv\pm\{7, 8, 9, 10, 11, 12, 13, 14, 15, 17, 19, 29\}\pmod {60},\; \text{and} \\
&T\equiv\pm\{1, 7, 10, 12, 15, 17, 19, 21, 23, 26, 28, 29\} \pmod{60}.\; \text{Then}\\
&p(S,n)=p(T,n-7),\; \text{for all}\; n\geq 7.
\end{align*}
\end{theorem}
\begin{proof}
We prove (i) which is equivalent to
\begin{align}
&[7,9,26,28,29:60]-q[8,14,17,19,21:60]\notag\\
&=[1, 7, 8, 9, 10, 11, 12, 13, 14, 15, 17, 19, 21, 23, 26, 28,
29:60].\label{sn1}
\end{align}
Similarly by the first part of Theorem \eqref{60.5}, we have
\begin{align}
&[6,9,21,28:60]-q[8,11,14,19:60]\notag\\&=[1, 6, 7, 8, 9, 10, 11,
13, 14, 17, 19, 20, 21, 23, 28, 29:60].\label{sn2}
\end{align}
By \eqref{sn2}, the equation \eqref{sn1} is equivalent to
\begin{align}
&[6,20:60]\big\{[7,9,26,28,29:60]-q[8,14,17,19,21:60]\big\}\notag\\
&=[12,15,26:60]\big\{[6,9,21,28:60]-q[8,11,14,19:60]\big\}.\label{sn3}
\end{align}
After regrouping  the the terms in \eqref{sn3}, we  are led to
prove
\begin{align}
&[6,9,28:60]\big\{[7,20,26,29:60]-[12,15,21,26:60]\big\}\notag\\
&=q[8,14,19:60]\big\{[6,17,20,21:60]-[11,12,15,26:60]\big\}.\label{sn4}
\end{align}
Next, we apply \eqref{four} twice with the sets of parameters
$[1,6,13,21,27]$, $[1,6,12,18,21]$ and with $q$ replaced by $q^{60}$
in each instance to find that
\begin{align*}
& [7,20,26,29:60]-[12,15,21,26:60]=-q^7[5,8,14,19:60]\;\text{and}\\
&[6,17,20,21:60]-[11,12,15,26:60]=-q^6[5,6,9,28:60].
\end{align*}
These two identities readily imply \eqref{sn4}. Hence, the proof of
(i) is complete.

\end{proof}
\begin{theorem}\label{60.7}

\begin{align*}
\mathrm{(i)}\hspace{0.5in}&\text{Let}\; S\equiv\pm\{1, 4, 7, 10, 11, 13, 17, 18, 19, 24, 25, 27\}\pmod {60},\; \text{and} \\
&T\equiv\pm\{1, 3, 7, 10, 11, 16, 17, 18, 23, 24, 25, 29\} \pmod{60}.\; \text{Then}\\
&p(S,n)=p(T,n-1),\; \text{for all}\; n\geq 1.
\end{align*}
\begin{align*}
\mathrm{(ii)}\hspace{0.5in}&\text{Let}\; S\equiv\pm\{1, 5, 6, 7, 8, 10, 11, 12, 17, 19, 21, 23\}\pmod {60},\; \text{and} \\
&T\equiv\pm\{1, 5, 6, 7, 9, 10, 11, 12, 13, 17, 28, 29\}\pmod{60}.\; \text{Then}\\
&p(S,n)=p(T,n),\; \text{for all}\; n\neq 8.
\end{align*}
\begin{align*}
\mathrm{(iii)}\hspace{0.5in}&\text{Let}\; S\equiv\pm\{6, 7, 8, 9, 10, 11, 12, 13, 19, 23, 25, 29\}\pmod {60},\; \text{and} \\
&T\equiv\pm\{1, 6, 10, 12, 13, 17, 19, 21, 23, 25, 28, 29\} \pmod{60}.\; \text{Then}\\
&p(S,n)=p(T,n-6),\; \text{for all}\; n\geq 6.
\end{align*}
\begin{align*}
\mathrm{(iv)}\hspace{0.5in}&\text{Let}\; S\equiv\pm\{3, 4, 5, 10, 11, 13, 17, 18, 19, 23, 24, 29\}\pmod {60},\; \text{and} \\
&T\equiv\pm\{1, 5, 7, 10, 13, 16, 18, 19, 23, 24, 27, 29\} \pmod{60}.\; \text{Then}\\
&p(S,n)=p(T,n-3),\; \text{for all}\; n\geq 3.
\end{align*}
\end{theorem}
\begin{proof}
We prove (i) which is equivalent to
\begin{align}
&[3,16,23,29:60]-q[4,13,19,27:60]\notag\\
&=[1, 3, 4, 7, 10, 11, 13, 16, 17, 18, 19, 23, 24, 25, 27,
29:60].\label{snn1}
\end{align}
Similarly by the second part of Theorem \eqref{60.5}, we have
\begin{align}
&[3,16,18,27:60]-q[4,13,17,22:60]\notag\\&=[1, 3, 4, 7, 10, 11,
13, 16, 17, 18, 19, 20, 22, 23, 27, 29:60].\label{snn2}
\end{align}
Therefore, by \eqref{snn2}, it suffices to prove that
\begin{align}
&[20,22:60]\big\{[3,16,23,29:60]-q[4,13,19,27:60]\big\}\notag\\
&=[24,25:60]\big\{[3,16,18,27:60]-q[4,13,17,22:60]\big\}.\label{snn3}
\end{align}
After regrouping  the the terms in \eqref{snn3}, we  are led to
prove
\begin{align}
&[3,16:60]\big\{[20,22,23,29:60]-[18,24,25,27:60]\big\}\notag\\
&=-q[4,13:60]\big\{[17,22,24,25:60]-[19,20,22,27:60]\big\}.\label{snn4}
\end{align}
Next, we apply \eqref{four} twice with the sets of parameters
$[1,3,21,25,26]$, $[1,3,20,23,25]$ and with $q$ replaced by $q^{60}$
in each instance to find that
\begin{align*}
& [20,22,23,29:60]-[18,24,25,27:60]=q^{18}[2,4,5,13:60]\;\text{and}\\
&[17,22,24,25:60]-[19,20,22,27:60]=-q^{17}[2,3,5,16:60].
\end{align*}
From these two identities  \eqref{snn4} trivially follows and so
the proof of (i) is complete.
\end{proof}
\begin{theorem}\label{60.8}

\begin{align*}
\mathrm{(i)}\hspace{0.5in}&\text{Let}\; S\equiv\pm\{1, 6, 7, 8, 11, 13, 18, 19, 20, 21, 23, 25\}\pmod {60},\; \text{and} \\
&T\equiv\pm\{1, 5, 6, 7, 13, 17, 18, 19, 20, 21, 28, 29\} \pmod{60}.\; \text{Then}\\
&p(S,n)=p(T,n-1),\; \text{for all}\; n\geq 1.
\end{align*}
\begin{align*}
\mathrm{(ii)}\hspace{0.5in}&\text{Let}\; S\equiv\pm\{4, 5, 6, 7, 11, 13, 17, 18, 19, 20, 27, 29\}\pmod {60},\; \text{and} \\
&T\equiv\pm\{1, 6, 7, 11, 13, 16, 18, 20, 23, 25, 27, 29\}\pmod{60}.\; \text{Then}\\
&p(S,n)=p(T,n-4),\; \text{for all}\; n\geq 4.
\end{align*}
\begin{align*}
\mathrm{(iii)}\hspace{0.5in}&\text{Let}\; S\equiv\pm\{5, 6, 7, 8, 9, 11, 17, 18, 19, 20, 23, 29\}\pmod {60},\; \text{and} \\
&T\equiv\pm\{1, 6, 9, 11, 13, 17, 18, 20, 23, 25, 28, 29\} \pmod{60}.\; \text{Then}\\
&p(S,n)=p(T,n-5),\; \text{for all}\; n\geq 5.
\end{align*}
\begin{align*}
\mathrm{(iv)}\hspace{0.5in}&\text{Let}\; S\equiv\pm\{1, 3, 4, 6, 11, 13, 17, 18, 19, 20, 23, 25\}\pmod {60},\; \text{and} \\
&T\equiv\pm\{1, 3, 5, 6, 7, 16, 17, 18, 19, 20, 23, 29\} \pmod{60}.\; \text{Then}\\
&p(S,n)=p(T,n),\; \text{for all}\; n\neq 4.
\end{align*}
\end{theorem}
\begin{proof}
We prove (i) which is equivalent to
\begin{align}
&[5,17,28,29:60]-q[8,11,23,25:60]\notag\\
&=[1, 5, 6, 7, 8, 11, 13, 17, 18, 19, 20, 21, 23, 25, 28,
29:60].\label{sv1}
\end{align}
Similarly by the second part of Theorem \eqref{60.3}, we have
\begin{align}
&[5,16,25:60]-q[8,11,19:60]\notag\\&=[1, 5, 6, 7, 8, 11, 13, 16,
17, 18, 19, 23, 24, 25, 29:60].\label{sv2}
\end{align}
Therefore, we need to prove that
\begin{align}
&[20,21,28:60]\big\{[5,16,25:60]-q[8,11,19:60]\big\}\notag\\
&=[16,24:60]\big\{[5,17,28,29:60]-q[8,11,23,25:60]\big\}.\label{sv3}
\end{align}
After regrouping  the the terms in \eqref{sv3}, we  are led to
prove
\begin{align}
&[5,16:60]\big\{[20,21,25,28:60]-[17,24,28,29:60]\big\}\notag\\
&=-q[8,11:60]\big\{[16,23,24,25:60]-[19,20,21,28:60]\big\}.\label{sv4}
\end{align}
Next, we apply \eqref{four} twice with the sets of parameters
$[1,4,21,25,29]$, $[1,4,20,24,25]$ and with $q$ replaced by $q^{60}$
in each instance to find that
\begin{align*}
& [20,21,25,28:60]-[17,24,28,29:60]=q^{17}[3,4,8,11:60]\;\text{and}\\
&[16,23,24,25:60]-[19,20,21,28:60]=-q^{16}[3,4,5,16:60].
\end{align*}
From these two identities  \eqref{sv4} trivially follows and so
the proof of (i) is complete.
\end{proof}

\section{Modulus $M = 62$}
There is only one class  of identities under multiplication
by the group $U(62)$. \\
\begin{theorem}\label{62.1}

\begin{align*}
\mathrm{(i)}\hspace{0.5in}&\text{Let}\; S\equiv\pm\{1, 6, 7, 11, 12, 15, 16, 17, 19, 20, 25, 29\}\pmod {62},\; \text{and} \\
&T\equiv\pm\{1, 5, 6, 11, 13, 14, 17, 18, 21, 25, 29, 30\} \pmod{62}.\; \text{Then}\\
&p(S,n)=p(T,n-1),\; \text{for all}\; n\geq 1.
\end{align*}
\begin{align*}
\mathrm{(ii)}\hspace{0.5in}&\text{Let}\; S\equiv\pm\{2, 3, 5, 11, 13, 14, 17, 18, 21, 25, 26, 29\}\pmod {62},\; \text{and} \\
&T\equiv\pm\{1, 3, 8, 11, 13, 15, 18, 20, 23, 25, 28, 29\} \pmod{62}.\; \text{Then}\\
&p(S,n)=p(T,n-2),\; \text{for all}\; n\geq 2.
\end{align*}
\begin{align*}
\mathrm{(iii)}\hspace{0.5in}&\text{Let}\; S\equiv\pm\{1, 3, 5, 7, 8, 19, 21, 23, 25, 26, 28, 30\}\pmod {62},\; \text{and} \\
&T\equiv\pm\{1, 2, 5, 7, 13, 18, 21, 23, 24, 27, 29, 30\} \pmod{62}.\; \text{Then}\\
&p(S,n)=p(T,n-1),\; \text{for all}\; n\geq 1.
\end{align*}
\begin{align*}
\mathrm{(iv)}\hspace{0.5in}&\text{Let}\; S\equiv\pm\{5, 7, 9, 11, 12, 13, 15, 16, 17, 19, 20, 22\}\pmod {62},\; \text{and} \\
&T\equiv\pm\{2, 5, 7, 11, 15, 17, 20, 23, 24, 26, 27, 29\} \pmod{62}.\; \text{Then}\\
&p(S,n)=p(T,n-5),\; \text{for all}\; n\geq 5.
\end{align*}
\begin{align*}
\mathrm{(v)}\hspace{0.5in}&\text{Let}\; S\equiv\pm\{2, 3, 7, 8, 9, 13, 17, 22, 23, 24, 25, 29\}\pmod {62},\; \text{and} \\
&T\equiv\pm\{1, 6, 8, 9, 11, 13, 15, 16, 20, 23, 25, 29\} \pmod{62}.\; \text{Then}\\
&p(S,n)=p(T,n),\; \text{for all}\; n\neq 1.
\end{align*}
\begin{align*}
\mathrm{(vi)}\hspace{0.5in}&\text{Let}\; S\equiv\pm\{1, 3, 4, 7, 9, 11, 12, 17, 19, 20, 27, 30\}\pmod {62},\; \text{and} \\
&T\equiv\pm\{1, 3, 4, 8, 9, 10, 11, 15, 21, 23, 27, 28\} \pmod{62}.\; \text{Then}\\
&p(S,n)=p(T,n),\; \text{for all}\; n\neq 7.
\end{align*}
\begin{align*}
\mathrm{(vii)}\hspace{0.5in}&\text{Let}\; S\equiv\pm\{3, 4, 5, 13, 14, 15, 16, 17, 18, 19, 25, 27\}\pmod {62},\; \text{and} \\
&T\equiv\pm\{1, 5, 9, 12, 13, 15, 16, 19, 22, 27, 29, 30\} \pmod{62}.\; \text{Then}\\
&p(S,n)=p(T,n-3),\; \text{for all}\; n\geq 3.
\end{align*}
\begin{align*}
\mathrm{(viii)}\hspace{0.5in}&\text{Let}\; S\equiv\pm\{1, 3, 5, 7, 9, 13, 15, 16, 21, 22, 24, 28\}\pmod {62},\; \text{and} \\
&T\equiv\pm\{1, 3, 6, 7, 8, 10, 15, 19, 21, 23, 25, 28\} \pmod{62}.\; \text{Then}\\
&p(S,n)=p(T,n),\; \text{for all}\; n\neq 5.
\end{align*}
\begin{align*}
\mathrm{(ix)}\hspace{0.5in}&\text{Let}\; S\equiv\pm\{1, 3, 4, 9, 10, 14, 15, 17, 21, 22, 23, 27\}\pmod {62},\; \text{and} \\
&T\equiv\pm\{1, 3, 5, 7, 9, 13, 17, 18, 21, 22, 24, 30\}\pmod{62}.\; \text{Then}\\
&p(S,n)=p(T,n),\; \text{for all}\; n\neq 4.
\end{align*}
\begin{align*}
\mathrm{(x)}\hspace{0.5in}&\text{Let}\; S\equiv\pm\{6, 7, 8, 9, 10, 11, 13, 19, 20, 21, 23, 25\}\pmod {62},\; \text{and} \\
&T\equiv\pm\{1, 7, 10, 12, 13, 18, 19, 21, 23, 27, 29, 30\}\pmod{62}.\; \text{Then}\\
&p(S,n)=p(T,n-6),\; \text{for all}\; n\geq 6.
\end{align*}
\begin{align*}
\mathrm{(xi)}\hspace{0.5in}&\text{Let}\; S\equiv\pm\{2, 6, 7, 10, 11, 15, 16, 17, 19, 21, 25, 29\}\pmod {62},\; \text{and} \\
&T\equiv\pm\{2, 4, 5, 11, 14, 15, 17, 21, 23, 26, 27, 29\}\pmod{62}.\; \text{Then}\\
&p(S,n)=p(T,n-2),\; \text{for all}\; n\geq 2.
\end{align*}
\begin{align*}
\mathrm{(xii)}\hspace{0.5in}&\text{Let}\; S\equiv\pm\{5, 8, 9, 11, 12, 13, 14, 15, 17, 19, 20, 23\}\pmod {62},\; \text{and} \\
&T\equiv\pm\{3, 4, 5, 14, 15, 17, 19, 23, 25, 26, 27, 28\}\pmod{62}.\; \text{Then}\\
&p(S,n)=p(T,n-5),\; \text{for all}\; n\geq 5.
\end{align*}
\begin{align*}
\mathrm{(xiii)}\hspace{0.5in}&\text{Let}\; S\equiv\pm\{3, 4, 5, 9, 10, 11, 19, 21, 25, 26, 27, 28\}\pmod {62},\; \text{and} \\
&T\equiv\pm\{1, 5, 6, 9, 15, 16, 19, 22, 25, 26, 27, 29\}\pmod{62}.\; \text{Then}\\
&p(S,n)=p(T,n-3),\; \text{for all}\; n\geq 3.
\end{align*}
\begin{align*}
\mathrm{(xiv)}\hspace{0.5in}&\text{Let}\; S\equiv\pm\{4, 6, 7, 9, 10, 11, 13, 21, 23, 24, 25, 27\}\pmod {62},\; \text{and} \\
&T\equiv\pm\{2, 3, 7, 13, 14, 17, 18, 23, 24, 25, 27, 29\}\pmod{62}.\; \text{Then}\\
&p(S,n)=p(T,n-4),\; \text{for all}\; n\geq 4.
\end{align*}
\begin{align*}
\mathrm{(xv)}\hspace{0.5in}&\text{Let}\; S\equiv\pm\{2, 3, 5, 9, 11, 12, 19, 21, 26, 27, 28, 29\}\pmod {62},\; \text{and} \\
&T\equiv\pm\{1, 3, 7, 9, 12, 17, 19, 22, 24, 27, 29, 30\}\pmod{62}.\; \text{Then}\\
&p(S,n)=p(T,n-2),\; \text{for all}\; n\geq 2.
\end{align*}
\end{theorem}
\begin{proof}
The identities (i)--(xv) follow from \eqref{four2} with the
following sets of parameters
\begin{align*}
&[1, 3, 2, 8, 13],\,[1, 3, 4, 17, 24],\,[1, 2, 4, 14, 27],\,[1, 6,
8, 24, 25],\, [1, 3, 4, 10, 12],\,\\ &[1, 2, 9, 12, 13],\,[1, 5,
2, 7, 14],\,[1, 2, 7, 10, 14],\,[1, 2, 6, 9, 15],\,[1, 8, 2, 10,
12],\\ &[1, 3, 7, 22, 24], [1, 6, 9, 24, 26],[1, 5, 2, 7, 11],[1,
7, 3, 10, 11], [1, 4, 2,6, 11],
\end{align*}
and with $q$ replaced by $q^{31}$ in each instance.
\end{proof}
\section{Modulus $M = 64$}
There is only one class  of identities under multiplication by the
group $U(64)$. \\
\begin{theorem}\label{64.1}

\begin{align*}
\mathrm{(i)}\hspace{0.5in}&\text{Let}\; S\equiv\pm\{4, 5, 6, 7, 11, 13, 19, 21, 22, 25, 27, 28\}\pmod {64},\; \text{and} \\
&T\equiv\pm\{1, 6, 7, 11, 15, 17, 20, 21, 25, 26, 28, 31\}\pmod{64}.\; \text{Then}\\
&p(S,n)=p(T,n-4),\; \text{for all}\; n\geq 4.
\end{align*}
\begin{align*}
\mathrm{(ii)}\hspace{0.5in}&\text{Let}\; S\equiv\pm\{1, 2, 7, 11, 12, 15, 17, 18, 20, 21, 25, 31\}\pmod {64},\; \text{and} \\
&T\equiv\pm\{1, 3, 4, 11, 13, 14, 18, 19, 20, 21, 29, 31\}\pmod{64}.\; \text{Then}\\
&p(S,n)=p(T,n),\; \text{for all}\; n\neq 2.
\end{align*}
\begin{align*}
\mathrm{(iii)}\hspace{0.5in}&\text{Let}\; S\equiv\pm\{2, 3, 5, 9, 11, 12, 21, 23, 27, 28, 29, 30\}\pmod {64},\; \text{and} \\
&T\equiv\pm\{1, 3, 7, 9, 12, 18, 20, 23, 25, 29, 30, 31\}\pmod{64}.\; \text{Then}\\
&p(S,n)=p(T,n-2),\; \text{for all}\; n\geq 2.
\end{align*}
\begin{align*}
\mathrm{(iv)}\hspace{0.5in}&\text{Let}\; S\equiv\pm\{4, 7, 9, 10, 12, 13, 15, 17, 19, 22, 23, 25\}\pmod {64},\; \text{and} \\
&T\equiv\pm\{3, 4, 5, 13, 15, 17, 19, 22, 26, 27, 28, 29\}\pmod{64}.\; \text{Then}\\
&p(S,n)=p(T,n-4),\; \text{for all}\; n\geq 4.
\end{align*}
\begin{align*}
\mathrm{(v)}\hspace{0.5in}&\text{Let}\; S\equiv\pm\{1, 3, 4, 6, 10, 11, 13, 19, 21, 28, 29, 31\}\pmod {64},\; \text{and} \\
&T\equiv\pm\{1, 3, 4, 7, 9, 10, 12, 22, 23, 25, 29, 31\}\pmod{64}.\; \text{Then}\\
&p(S,n)=p(T,n),\; \text{for all}\; n\neq 6.
\end{align*}

\begin{align*}
\mathrm{(vi)}\hspace{0.5in}&\text{Let}\; S\equiv\pm\{2, 7, 9, 12, 13, 14, 15, 17, 19, 20, 23, 25\}\pmod {64},\; \text{and} \\
&T\equiv\pm\{2, 5, 7, 11, 12, 13, 19, 21, 25, 27, 28, 30\}\pmod{64}.\; \text{Then}\\
&p(S,n)=p(T,n-2),\; \text{for all}\; n\geq 2.
\end{align*}
\begin{align*}
\mathrm{(vii)}\hspace{0.5in}&\text{Let}\; S\equiv\pm\{3, 4, 5, 13, 14, 15, 17, 18, 19, 20, 27, 29\}\pmod {64},\; \text{and} \\
&T\equiv\pm\{1, 5, 9, 12, 14, 15, 17, 20, 23, 27, 30, 31\}\pmod{64}.\; \text{Then}\\
&p(S,n)=p(T,n-3),\; \text{for all}\; n\geq 3.
\end{align*}
\begin{align*}
\mathrm{(viii)}\hspace{0.5in}&\text{Let}\; S\equiv\pm\{3, 4, 5, 9, 10, 11, 21, 23, 26, 27, 28, 29\}\pmod {64},\; \text{and} \\
&T\equiv\pm\{1, 5, 6, 9, 15, 17, 20, 23, 26, 27, 28, 31\}\pmod{64}.\; \text{Then}\\
&p(S,n)=p(T,n-3),\; \text{for all}\; n\geq 3.
\end{align*}

\end{theorem}
\begin{proof}
The identities (i)--(viii) follow from \eqref{four2} with the
following sets of parameters
\begin{align*}
&[1, 5, 6, 18, 27],\;[1, 2, 4, 8, 15],\;[1, 4, 2, 6, 11],\; [1, 5,
8, 23, 27],\;[1, 2, 8, 11, 12], \\ &[1, 5, 3, 10, 15],\;  [1, 5, 2,
7, 14], [1, 5, 2, 7, 11].
\end{align*}
and with $q$ replaced by $q^{32}$ in each instance.
\end{proof}

\section{Modulus $M = 66$}
There are two classes  of identities under multiplication by the
group $U(66)$.\\
\begin{theorem}\label{66.1}

\begin{align*}
\mathrm{(i)}\hspace{0.5in}&\text{Let}\; S\equiv\pm\{7, 8, 9, 10, 11, 12, 13, 19, 21, 22, 23, 25\}\pmod {66},\; \text{and} \\
&T\equiv\pm\{1, 9, 11, 12, 14, 19, 21, 22, 25, 29, 31, 32\} \pmod{66}.\; \text{Then}\\
&p(S,n)=p(T,n-7),\; \text{for all}\; n\geq 7.
\end{align*}
\begin{align*}
\mathrm{(ii)}\hspace{0.5in}&\text{Let}\; S\equiv\pm\{4, 5, 6, 7, 11, 13, 21, 22, 23, 27, 28, 29\}\pmod {66},\; \text{and} \\
&T\equiv\pm\{1, 6, 7, 11, 16, 17, 21, 22, 26, 27, 29, 31\} \pmod{66}.\; \text{Then}\\
&p(S,n)=p(T,n-4),\; \text{for all}\; n\geq 4.
\end{align*}
\begin{align*}
\mathrm{(iii)}\hspace{0.5in}&\text{Let}\; S\equiv\pm\{1, 3, 5, 7, 11, 15, 18, 19, 22, 23, 26, 32\}\pmod {66},\; \text{and} \\
&T\equiv\pm\{1, 3, 4, 10, 11, 15, 17, 18, 22, 23, 25, 29\} \pmod{66}.\; \text{Then}\\
&p(S,n)=p(T,n),\; \text{for all}\; n\neq 4.
\end{align*}
\begin{align*}
\mathrm{(iv)}\hspace{0.5in}&\text{Let}\; S\equiv\pm\{5, 7, 9, 11, 13, 15, 16, 17, 19, 20, 22, 24\}\pmod {66},\; \text{and} \\
&T\equiv\pm\{2, 5, 9, 11, 15, 17, 22, 24, 25, 28, 29, 31\} \pmod{66}.\; \text{Then}\\
&p(S,n)=p(T,n-5),\; \text{for all}\; n\geq 5.
\end{align*}
\begin{align*}
\mathrm{(v)}\hspace{0.5in}&\text{Let}\; S\equiv\pm\{2, 3, 5, 11, 13, 14, 19, 22, 23, 27, 30, 31\}\pmod {66},\; \text{and} \\
&T\equiv\pm\{1, 3, 8, 11, 13, 17, 20, 22, 25, 27, 30, 31\} \pmod{66}.\; \text{Then}\\
&p(S,n)=p(T,n-2),\; \text{for all}\; n\geq 2.
\end{align*}
\end{theorem}
\begin{proof}
The identities (i)--(v) follow from \eqref{four2} with the following
sets of parameters
\begin{equation*}
[1, 9, 2, 11, 13],\, [1, 5, 6, 18, 28],\, [1, 2, 6, 9, 16],\,[1,
6, 8, 25, 26],\, [1, 3, 4, 17, 26],
\end{equation*}
and with $q$ replaced by $q^{33}$ in each instance.
\end{proof}
\begin{theorem}\label{66.2}

\begin{align*}
\mathrm{(i)}\hspace{0.5in}&\text{Let}\; S\equiv\pm\{1, 10, 11, 12, 13, 14, 15, 16, 19, 21, 23, 25\}\pmod {66},\; \text{and} \\
&T\equiv\pm\{1, 9, 10, 11, 12, 13, 14, 15, 25, 29, 31, 32\} \pmod{66}.\; \text{Then}\\
&p(S,n)=p(T,n-1),\; \text{for all}\; n\geq 1.
\end{align*}
\begin{align*}
\mathrm{(ii)}\hspace{0.5in}&\text{Let}\; S\equiv\pm\{1, 4, 5, 6, 7, 9, 11, 13, 16, 21, 23, 28\}\pmod {66},\; \text{and} \\
&T\equiv\pm\{1, 4, 5, 6, 7, 9, 11, 14, 16, 17, 27, 29\} \pmod{66}.\; \text{Then}\\
&p(S,n)=p(T,n),\; \text{for all}\; n\neq 13.
\end{align*}
\begin{align*}
\mathrm{(iii)}\hspace{0.5in}&\text{Let}\; S\equiv\pm\{3, 4, 5, 7, 11, 18, 19, 23, 25, 26, 27, 32\}\pmod {66},\; \text{and} \\
&T\equiv\pm\{1, 4, 7, 11, 15, 18, 20, 23, 25, 27, 29, 32\} \pmod{66}.\; \text{Then}\\
&p(S,n)=p(T,n-3),\; \text{for all}\; n\geq 3.
\end{align*}
\begin{align*}
\mathrm{(iv)}\hspace{0.5in}&\text{Let}\; S\equiv\pm\{2, 3, 5, 7, 11, 13, 15, 16, 19, 20, 24, 29\}\pmod {66},\; \text{and} \\
&T\equiv\pm\{2, 3, 5, 9, 10, 11, 13, 16, 17, 24, 29, 31\} \pmod{66}.\; \text{Then}\\
&p(S,n)=p(T,n),\; \text{for all}\; n\neq 7.
\end{align*}
\begin{align*}
\mathrm{(v)}\hspace{0.5in}&\text{Let}\; S\equiv\pm\{5, 6, 7, 8, 9, 11, 17, 23, 26, 27, 28, 29\}\pmod {66},\; \text{and} \\
&T\equiv\pm\{1, 6, 9, 11, 16, 17, 21, 23, 26, 28, 29, 31\} \pmod{66}.\; \text{Then}\\
&p(S,n)=p(T,n-5),\; \text{for all}\; n\geq 5.
\end{align*}
\begin{align*}
\mathrm{(vi)}\hspace{0.5in}&\text{Let}\; S\equiv\pm\{2, 5, 8, 11, 13, 14, 17, 19, 21, 23, 27, 30\}\pmod {66},\; \text{and} \\
&T\equiv\pm\{2, 3, 8, 11, 13, 17, 19, 21, 25, 26, 30, 31\} \pmod{66}.\; \text{Then}\\
&p(S,n)=p(T,n-2),\; \text{for all}\; n\geq 2.
\end{align*}
\begin{align*}
\mathrm{(vii)}\hspace{0.5in}&\text{Let}\; S\equiv\pm\{7, 8, 9, 10, 11, 12, 13, 15, 19, 23, 31, 32\}\pmod {66},\; \text{and} \\
&T\equiv\pm\{1, 8, 11, 12, 15, 19, 21, 23, 25, 28, 31, 32\}\pmod{66}.\; \text{Then}\\
&p(S,n)=p(T,n-7),\; \text{for all}\; n\geq 7.
\end{align*}
\begin{align*}
\mathrm{(viii)}\hspace{0.5in}&\text{Let}\; S\equiv\pm\{3, 4, 5, 11, 13, 14, 19, 20, 21, 25, 30, 31\}\pmod {66},\; \text{and} \\
&T\equiv\pm\{1, 5, 8, 11, 14, 17, 20, 21, 25, 27, 30, 31\}\pmod{66}.\; \text{Then}\\
&p(S,n)=p(T,n-3),\; \text{for all}\; n\geq 3.
\end{align*}
\begin{align*}
\mathrm{(ix)}\hspace{0.5in}&\text{Let}\; S\equiv\pm\{1, 3, 4, 10, 11, 17, 18, 19, 25, 26, 27, 29\}\pmod {66},\; \text{and} \\
&T\equiv\pm\{1, 2, 7, 10, 11, 15, 18, 19, 23, 26, 27, 29\}\pmod{66}.\; \text{Then}\\
&p(S,n)=p(T,n-2),\; \text{for all}\; n\geq 2.
\end{align*}
\begin{align*}
\mathrm{(x)}\hspace{0.5in}&\text{Let}\; S\equiv\pm\{3, 5, 7, 9, 11, 13, 17, 20, 24, 28, 31, 32\}\pmod {66},\; \text{and} \\
&T\equiv\pm\{2, 3, 7, 11, 15, 17, 20, 24, 25, 28, 29, 31\}\pmod{66}.\; \text{Then}\\
&p(S,n)=p(T,n-3),\; \text{for all}\; n\geq 3.
\end{align*}
\end{theorem}
\begin{proof}
We prove (i) which is equivalent to
\begin{align}
&[9,29,31,32:66]-q[16,19,21,23:66]\notag\\
&=[1, 9, 10, 11, 12, 13, 14, 15, 16, 19, 21, 23, 25, 29, 31,
32:66].\label{sb1}
\end{align}
Our proof uses the first part of Theorem \eqref{66.1} that is
equivalent to
\begin{align}
&[1,14,29,31,32:66]-q^7[7,8,10,13,23:66]\notag\\&=[1, 7, 8, 9, 10,
11, 12, 13, 14, 19, 21, 22, 23, 25, 29, 31, 32:66].\label{sb2}
\end{align}
Therefore, it suffices to prove that
\begin{align}
&[15,16:66]\big\{[1,14,29,31,32:66]-q^7[7,8,10,13,23:66]\big\}\notag\\
&=[7,8,22:66]\big\{[9,29,31,32:66]-q[16,19,21,23:66]\big\}.\label{sy3}
\end{align}
After regrouping  the the terms in \eqref{sy3}, we  are led to
prove
\begin{align}
&[29,31,32:66]\big\{[1,14,15,16:66]-[7,8,9,22:66]\big\}\notag\\
&=-q[7,8,23:66]\big\{[16,19,21,22:66]-q^6[10,13,15,16:66]\big\}.\label{sy4}
\end{align}
Next, we apply \eqref{four} twice with the sets of parameters
$[1,2,24,8,9]$, $[1,11,17,30,32]$ and with $q$ replaced by $q^{66}$
in both instance to find that
\begin{align*}
&[1,14,15,16:66]-[7,8,9,22:66]=-q[6,7,8,23:66]\;\text{and}\\
&[16,19,21,22:66]-q^6[10,13,15,16:66]=[6,29,31,32:66].
\end{align*}
From these two identities  \eqref{sy4} readily follows and so the
proof of (i) is complete.
\end{proof}
\section{Modulus $M = 68$}
There is only one class  of identities under multiplication by the
group $U(68)$. \\
\begin{theorem}\label{68.1}

\begin{align*}
\mathrm{(i)}\hspace{0.5in}&\text{Let}\; S\equiv\pm\{1, 7, 8, 12, 13, 15, 18, 19, 21, 22, 27, 33\}\pmod {68},\; \text{and} \\
&T\equiv\pm\{1, 6, 7, 11, 14, 15, 19, 20, 23, 27, 32, 33\}\pmod{68}.\; \text{Then}\\
&p(S,n)=p(T,n-1),\; \text{for all}\; n\geq 1.
\end{align*}
\begin{align*}
\mathrm{(ii)}\hspace{0.5in}&\text{Let}\; S\equiv\pm\{2, 3, 5, 11, 13, 14, 21, 23, 24, 29, 31, 32\}\pmod {68},\; \text{and} \\
&T\equiv\pm\{1, 3, 8, 11, 13, 18, 21, 23, 26, 28, 31, 33\}\pmod{68}.\; \text{Then}\\
&p(S,n)=p(T,n-2),\; \text{for all}\; n\geq 2.
\end{align*}
\begin{align*}
\mathrm{(iii)}\hspace{0.5in}&\text{Let}\; S\equiv\pm\{1, 3, 5, 7, 8, 22, 26, 27, 28, 29, 31, 33\}\pmod {68},\; \text{and} \\
&T\equiv\pm\{1, 2, 5, 7, 13, 21, 24, 27, 29, 30, 32, 33\}\pmod{68}.\; \text{Then}\\
&p(S,n)=p(T,n-1),\; \text{for all}\; n\geq 1.
\end{align*}
\begin{align*}
\mathrm{(iv)}\hspace{0.5in}&\text{Let}\; S\equiv\pm\{3, 7, 10, 11, 12, 15, 16, 18, 19, 23, 27, 31\}\pmod {68},\; \text{and} \\
&T\equiv\pm\{3, 4, 7, 9, 15, 19, 20, 25, 26, 27, 30, 31\}\pmod{68}.\; \text{Then}\\
&p(S,n)=p(T,n-3),\; \text{for all}\; n\geq 3.
\end{align*}
\begin{align*}
\mathrm{(v)}\hspace{0.5in}&\text{Let}\; S\equiv\pm\{1, 4, 5, 6, 9, 15, 19, 25, 26, 28, 29, 33\}\pmod {68},\; \text{and} \\
&T\equiv\pm\{1, 3, 5, 9, 10, 14, 16, 24, 25, 29, 31, 33\}\pmod{68}.\; \text{Then}\\
&p(S,n)=p(T,n),\; \text{for all}\; n\neq 3.
\end{align*}
\begin{align*}
\mathrm{(vi)}\hspace{0.5in}&\text{Let}\; S\equiv\pm\{4, 5, 6, 7, 9, 11, 20, 23, 25, 27, 29, 30\}\pmod {68},\; \text{and} \\
&T\equiv\pm\{2, 5, 9, 11, 12, 15, 16, 18, 19, 23, 25, 29\}\pmod{68}.\; \text{Then}\\
&p(S,n)=p(T,n),\; \text{for all}\; n\neq 2.
\end{align*}
\begin{align*}
\mathrm{(vii)}\hspace{0.5in}&\text{Let}\; S\equiv\pm\{7, 8, 9, 10, 11, 12, 13, 21, 22, 23, 25, 27\}\pmod {68},\; \text{and} \\
&T\equiv\pm\{1, 9, 11, 13, 14, 20, 21, 23, 25, 30, 32, 33\}\pmod{68}.\; \text{Then}\\
&p(S,n)=p(T,n-7),\; \text{for all}\; n\geq 7.
\end{align*}
\begin{align*}
\mathrm{(viii)}\hspace{0.5in}&\text{Let}\; S\equiv\pm\{2, 3, 9, 10, 13, 15, 16, 19, 21, 24, 25, 31\}\pmod {68},\; \text{and} \\
&T\equiv\pm\{3, 4, 5, 6, 13, 15, 19, 21, 22, 28, 29, 31\}\pmod{68}.\; \text{Then}\\
&p(S,n)=p(T,n),\; \text{for all}\; n\neq 2.
\end{align*}

\end{theorem}
\begin{proof}
The identities (i)--(viii) follow from \eqref{four2} with the
following sets of parameters
\begin{align*}
&[1, 3, 2, 9, 14],\;[1, 3, 4, 17, 27],\;[1, 2, 4, 14, 30],\;[1, 4,
8, 26, 27],\; [1, 2, 5, 10, 11],\\ &[1, 5, 7, 12, 16],\;[1, 9, 2,
11, 13],\;    [1, 4, 6, 10, 16]
\end{align*}
and with $q$ replaced by $q^{34}$ in each instance.
\end{proof}

\section{Modulus $M = 70$}
There is one class  of identities under multiplication by the group
$U(70)$. The identity given  in Theorem \ref{70.1} below is
equivalent to \eqref{rr2}.
\begin{theorem}\label{70.1}
\begin{align*}
&\text{Let}\\
&S\equiv\pm\{1, 3, 4, 5, 6, 7, 9, 11, 13, 14, 15, 16, 17, 19, 23, 24, 25,
26, 27, 28, 29, 31, 33, 34\}\pmod {70},\\
&\text{and} \\
&T\equiv\pm\{1, 2, 3, 5, 8, 9, 11, 12, 13, 14, 15, 17, 18, 19, 21, 22, 23, 25, 27, 28, 29, 31, 32, 33\} \pmod{70}.\\
&\text{Then}\qquad p(S,n)=p(T,n-1),\; \text{for all}\; n \geq 1.
\end{align*}
\end{theorem}

\section{Modulus $M = 72$}
There are three classes  of identities under multiplication by the
group $U(72)$.
\begin{theorem}\label{72.1}
\begin{align*}
\mathrm{(i)}\hspace{0.5in}&\text{Let}\; S\equiv\pm\{1, 3, 5, 7, 8, 26, 28, 29, 30, 31, 33, 35\}\pmod {72},\; \text{and} \\
&T\equiv\pm\{1, 2, 5, 7, 13, 23, 28, 29, 31, 32, 34, 35\} \pmod{72}.\; \text{Then}\\
&p(S,n)=p(T,n-1),\; \text{for all}\; n\geq 1.
\end{align*}
\begin{align*}
\mathrm{(ii)}\hspace{0.5in}&\text{Let}\; S\equiv\pm\{1, 4, 5, 6, 11, 14, 15, 21, 25, 31, 32, 35\}\pmod {72},\; \text{and} \\
&T\equiv\pm\{1, 4, 5, 7, 10, 11, 16, 25, 26, 29, 31, 35\} \pmod{72}.\; \text{Then}\\
&p(S,n)=p(T,n),\; \text{for all}\; n\neq 6.
\end{align*}
\begin{align*}
\mathrm{(iii)}\hspace{0.5in}&\text{Let}\; S\equiv\pm\{1, 7, 8, 13, 14, 17, 19, 20, 22, 23, 29, 35\}\pmod {72},\; \text{and} \\
&T\equiv\pm\{1, 6, 7, 13, 15, 16, 20, 21, 23, 29, 34, 35\} \pmod{72}.\; \text{Then}\\
&p(S,n)=p(T,n-1),\; \text{for all}\; n\geq 1.
\end{align*}
\begin{align*}
\mathrm{(iv)}\hspace{0.5in}&\text{Let}\; S\equiv\pm\{2, 3, 5, 11, 16, 17, 19, 20, 25, 30, 31, 33\}\pmod {72},\; \text{and} \\
&T\equiv\pm\{1, 5, 8, 11, 14, 17, 19, 20, 22, 25, 31, 35\} \pmod{72}.\; \text{Then}\\
&p(S,n)=p(T,n),\; \text{for all}\; n\neq 1.
\end{align*}
\begin{align*}
\mathrm{(v)}\hspace{0.5in}&\text{Let}\; S\equiv\pm\{4, 7, 10, 11, 13, 16, 17, 19, 23, 25, 26, 29\}\pmod {72},\; \text{and} \\
&T\equiv\pm\{3, 4, 7, 13, 17, 19, 22, 23, 29, 30, 32, 33\}\pmod{72}.\; \text{Then}\\
&p(S,n)=p(T,n-4),\; \text{for all}\; n\geq 4.
\end{align*}
\begin{align*}
\mathrm{(vi)}\hspace{0.5in}&\text{Let}\; S\equiv\pm\{6, 8, 10, 11, 13, 15, 17, 19, 21, 23, 25, 28\}\pmod {72},\; \text{and} \\
&T\equiv\pm\{2, 5, 11, 13, 17, 19, 23, 25, 28, 31, 32, 34\}\pmod{72}.\; \text{Then}\\
&p(S,n)=p(T,n-6),\; \text{for all}\; n\geq 6.
\end{align*}
\end{theorem}
\begin{proof}
The identities (i)--(vi) follow from \eqref{four2} with the
following sets of parameters
\begin{equation*}
[1, 2, 4, 14, 32],\,[1, 2, 8, 12, 13],\,[1, 3, 2, 9, 15],\,[1, 3,
4, 9, 15],\,[1, 5, 8, 25, 30],\,  [1, 9, 3, 13, 14],
\end{equation*}
and with $q$ replaced by $q^{36}$ in each instance.

\end{proof}

\begin{theorem}\label{72.2}
\begin{align}
&\text{Let}\notag\\
&S\equiv\pm\{1, 3, 5, 7, 8, 9, 10, 11, 12, 13, 14, 16, 17,
18, 19, 23, 25, 27, 29, 31, 32, 33, 34, 35\}\pmod {72},\notag\\
&\text{and}\notag\\
&T\equiv\pm\{1, 2, 5, 7, 8, 9, 11, 12, 13, 15, 16, 17, 18, 19, 21, 22,
23, 25, 26, 27, 29, 31, 32, 35\}\pmod{72}.\notag\\
&\text{Then}\qquad p(S,n)=p(T,n-1),\; \text{for all}\; n\geq 1.\label{24prog}
\end{align}

\end{theorem}

\begin{proof}
Our proof relies heavily on the additive properties of theta
functions, and so we adopt the notation of \cite{III}. We first
recall  Ramanujan's definition for a general theta function. Let

\begin{equation}\label{gt}
f(a,b) := \sum_{n=-\i}^{\i}a^{n(n+1)/2}b^{n(n-1)/2}, \qquad |ab| <
1.
\end{equation}
The function $ f(a,b) $ satisfies the well-known Jacobi triple
product identity \cite[p.~35, Entry 19]{III}
\begin{equation}\label{JB}
f(a,b) = (-a;ab)_{\i}(-b;ab)_{\i}(ab;ab)_{\i}.
\end{equation}
Observe that
\begin{equation*}
[a:q]=\df{f(-a,-q/a)}{(q:q)_\i}.
\end{equation*}
By separating the sum into even and odd index terms, we obtain
from \eqref{gt} that \cite[p.~48, Entry 31]{III}
\begin{equation}\label{2d}
f(a,b)=f(a^3b,ab^3)+af(ab^{-1},b^3a^5).
\end{equation}

We are now ready to prove Theorem \ref{72.2}. First, we express
\eqref{24prog} in its equivalent form
\begin{align}
[2,15,21,22,26:72]-q[3,10,14,33,34:72]=\df{[1,2,\cdots,35:72]}{[4,6,20,24,28,30:72]}.\label{hwg1}
\end{align}
For the left hand side of \eqref{hwg1},
\begin{align}
L(q)&:=[2,15,21,22,26:72]-q[3,10,14,33,34:72]\notag\\
&=[2:24][15:36]-q[10:24][3:36]\notag\\
&=\df{1}{(q^{24};q^{24})_\i(q^{36};q^{36})_\i}\big\{f(-q^2,-q^{22})f(-q^{15},-q^{21})-qf(-q^{10},-q^{14})f(-q^3,-q^{33})\big\}.\label{hwg2}
\end{align}
By \eqref{2d},
\begin{equation}
f(q,-q^2)=f(-q^5,-q^7)+qf(-q,-q^{11})\;\; \text{and}\;\;
f(-q,q^2)=f(-q^5,-q^7)-qf(-q,-q^{11}).\label{2du}
\end{equation}
In \eqref{hwg2}, we use  \eqref{2du} twice with $q$ replaced by
$q^2$ and $q^3$, we find that
\begin{align}
4q^2L(q)&=\df{1}{(q^{24};q^{24})_\i(q^{36};q^{36})_\i}\Big\{(f(q^3,-q^6)+f(-q^3,q^6))(f(q^2,-q^4)-f(-q^2,q^4))\notag\\
&-(f(q^3,-q^6)-f(-q^3,q^6))(f(q^2,-q^4)+f(-q^2,q^4))\Big\}\notag\\
&=\df{2}{(q^{24};q^{24})_\i(q^{36};q^{36})_\i}\big\{f(-q^3,q^6)f(q^2,-q^4)-f(q^3,-q^6)f(-q^2,q^4)\big\}.\label{hwg3}
\end{align}
From \eqref{JB}, after some elementary product manipulations, we
find that
\begin{equation*}
f(-q,q^2)=\df{f(q^3,q^3)}{(-q;q^2)_\i}\;\; \text{and}\;\;
f(q,-q^2)=(-q,-q)_\i=\df{f(q,q)}{(-q;q^2)_\i}.
\end{equation*}

From these two equations we conclude that
\begin{equation}\label{hwg4}
2q^2L(q)=\df{1}{(q^{24};q^{24})_\i(q^{36};q^{36})_\i
(-q^2;q^4)_\i(-q^3;q^6)_\i}\big\{f(q^2,q^2)f(q^9,q^9)-f(q^3,q^3)f(q^6,q^6)\big\}.
\end{equation}
While for the right hand side of \eqref{hwg1}, we find that
\begin{align}
R(q)&:=\df{[1,2,\cdots,35:72]}{[4,6,20,24,28,30:72]}\notag\\
&=\df{(q;q)_\i}{(q^{36};q^{72})_\i(q^{72};q^{72})_\i[6:36][4:24][24:72]}\notag\\
&=\df{(q;q)_\i(q^{18};q^{36})_\i(q^{12};q^{24})_\i(q^{72};q^{72})_\i}{(q^{36};q^{36})_\i(q^{6};q^{12})_\i(q^{24};q^{24})_\i}\notag\\
&=\df{(q;q)_\i
f(q^6,q^{30})}{(q^{36};q^{36})_\i(q^{24};q^{24})_\i},\label{hwg5}
\end{align}
after several applications of \eqref{JB}. By \eqref{hwg3} and
\eqref{hwg5}, we see that \eqref{hwg1} is equivalent to
\begin{align}
&f(q^2,q^2)f(q^9,q^9)-f(q^3,q^3)f(q^6,q^6)\notag\\&=2q^3(-q^2;q^4)_\i
(-q^3;q^6)_\i (q;q)_\i f(q^6,q^{30})=
f(-q,-q^3)f(q^6,q^{30})(-q^3;q^6)_\i,\label{hwg6}
\end{align}
where in the last step \eqref{JB} is used.

Employing  \eqref{qpp} with $q$ replaced by $q^6$ and $x$ replaced
by $q$, we find after some algebra that
\begin{equation}\label{psiqp}
f(q^9,q^9)-qf(q^3,q^{15})=f(-q,-q^3)(-q^3;q^6)_\i.
\end{equation}
We need two more identities to establish \eqref{hwg6} namely,
\begin{equation}\label{t3d}
f(q,q)=f(q^9,q^9)+2qf(q^3,q^{15}) \;\;\text{and}
\end{equation}
\begin{equation}\label{terr}
f(q,q)f(q^2,q^2)=f(q^3,q^3)f(q^6,q^6)+2qf(q,q^5)f(q^2,q^{10}).
\end{equation}
The identity \eqref{t3d} is the 3-dissection of the classical theta
function $\varphi(q)=f(q,q)$ (see the corollary in
\cite[p.~49]{III}) and the identity \eqref{terr} can be obtained
from the Theorem on page 73 of \cite{III} with the following
choice of parameters \\
$\epsilon_1 = \epsilon_2 = 0$, $\alpha=\beta=1$, $a=b=q^2$,
$c=d=q$, and $m=3$.

We return to the left hand side of  \eqref{hwg6} and use
\eqref{terr} with $q$ replaced by $q^3$, \eqref{t3d} with $q$
replaced by $q^2$ and \eqref{psiqp}, we conclude that
\begin{align*}
&f(q^2,q^2)f(q^9,q^9)-f(q^3,q^3)f(q^6,q^6)\\
&=f(q^2,q^2)f(q^9,q^9)-f(q^9,q^9)f(q^{18},q^{18})-2q^3f(q^3,q^{15})f(q^6,q^{30})\\
&=f(q^9,q^9)\big\{f(q^2,q^2)-f(q^{18},q^{18})\big\}-2q^3f(q^3,q^{15})f(q^6,q^{30})\\
&=2q^3f(q^9,q^9)f(q^6,q^{30})-2q^3f(q^3,q^{15})f(q^6,q^{30})\\
&=2q^2f(q^6,q^{30})\big\{f(q^9,q^9)-qf(q^3,q^{15})\big\} \\
&=2q^2f(q^6,q^{30})f(-q,-q^3)(-q^3;q^6)_\i.
\end{align*}
This is \eqref{hwg6} and so  the proof of \eqref{hwg1} is
complete.
\end{proof}
\begin{theorem}\label{72.3}

\begin{align*}
\mathrm{(i)}\hspace{0.5in}&\text{Let}\; S\equiv\pm\{1, 5, 9, 12, 13, 14, 16, 23, 27, 31, 34, 35\}\pmod {72},\; \text{and} \\
&T\equiv\pm\{1, 4, 9, 11, 13, 14, 23, 25, 27, 30, 32, 35\} \pmod{72}.\; \text{Then}\\
&p(S,n)=p(T,n-1),\; \text{for all}\; n\geq 1.
\end{align*}
\begin{align*}
\mathrm{(ii)}\hspace{0.5in}&\text{Let}\; S\equiv\pm\{2, 5, 6, 7, 9, 16, 17, 19, 20, 27, 29, 31\}\pmod {72},\; \text{and} \\
&T\equiv\pm\{2, 5, 7, 8, 9, 11, 12, 25, 26, 27, 29, 31\} \pmod{72}.\; \text{Then}\\
&p(S,n)=p(T,n),\; \text{for all}\; n\neq 6.
\end{align*}
\begin{align*}
\mathrm{(iii)}\hspace{0.5in}&\text{Let}\; S\equiv\pm\{5, 6, 7, 8, 9, 17, 19, 26, 27, 28, 29, 31\}\pmod {72},\; \text{and} \\
&T\equiv\pm\{1, 7, 9, 12, 17, 19, 22, 26, 27, 29, 32, 35\} \pmod{72}.\; \text{Then}\\
&p(S,n)=p(T,n-5),\; \text{for all}\; n\geq 5.
\end{align*}
\begin{align*}
\mathrm{(iv)}\hspace{0.5in}&\text{Let}\; S\equiv\pm\{1, 9, 10, 11, 12, 14, 17, 19, 25, 27, 32, 35\}\pmod {72},\; \text{and} \\
&T\equiv\pm\{1, 8, 9, 10, 11, 13, 23, 25, 27, 28, 30, 35\} \pmod{72}.\; \text{Then}\\
&p(S,n)=p(T,n-1),\; \text{for all}\; n\geq 1.
\end{align*}
\begin{align*}
\mathrm{(v)}\hspace{0.5in}&\text{Let}\; S\equiv\pm\{7, 8, 9, 10, 11, 12, 13, 23, 25, 27, 29, 34\}\pmod {72},\; \text{and} \\
&T\equiv\pm\{1, 9, 11, 13, 16, 20, 23, 25, 27, 30, 34, 35\} \pmod{72}.\; \text{Then}\\
&p(S,n)=p(T,n-7),\; \text{for all}\; n\geq 7.
\end{align*}
\begin{align*}
\mathrm{(vi)}\hspace{0.5in}&\text{Let}\; S\equiv\pm\{4, 5, 6, 7, 9, 17, 19, 22, 27, 29, 31, 32\}\pmod {72},\; \text{and} \\
&T\equiv\pm\{2, 5, 9, 12, 13, 16, 17, 19, 22, 23, 27, 31\} \pmod{72}.\; \text{Then}\\
&p(S,n)=p(T,n),\; \text{for all}\; n\neq 2.
\end{align*}
\end{theorem}
\begin{proof}
We prove (i) which is equivalent to
\begin{align}
&[4,11,25,30,32:72]-q[5,12,16,31,34:72]\notag\\
&=[1,4,5,9,11,12,13,14,16,23,25,27,30,31,32,34,35:72].\label{bt}
\end{align}
From  \eqref{four2} with $q$ replaced by $q^{36}$ and
$a,\,b,\,c,\,x,\,y,$ replaced by
$q,\,q^2,\,q^6,\,q^{12},\,q^{15}$, respectively, we find that
\begin{align*}
&\df{[8,28,34:72]}{[1,4,6,9,11,12,14,17,19,24,25,27,30,32,35:72]}\notag\\
&\;-\df{[32:72]}{[1,5,6,9,12,13,16,23,24,27,30,31,35:72]}=q^4.
\end{align*}
After clearing denominators, we arrive at
\begin{align}
&[5,8,13,16,23,28,31,34:72]-[4,11,14,17,19,25,32,32:72]\notag\\&
=q^4[1,4,5,6,9,11,12,13,14,16,17,19,23,24,25,27,30,31,32,35:72].\label{bt1}
\end{align}
By \eqref{bt1}, we see that \eqref{bt} is equivalent to
\begin{align}
&[34:72]\big\{[5,8,13,16,23,28,31,34:72]-[4,11,14,17,19,25,32,32:72]\big\}\notag\\
&=q^4[6,17,19,24:72]\big\{[4,11,25,30,32:72]-q[5,12,16,31,34:72]\big\}.\label{bt2}
\end{align}
After rearrangement of the terms of \eqref{bt2}, we are led the
prove
\begin{align}
&[4,11,17,19,25,32:72]\big\{[14,32,34:72]+q^4[6,24,30:72]\big\}\notag\\
&=[5,16,31,34:72]\big\{[8,13,23,28,34:72]+q^5[6,12,17,19,24:72]\big\}.\label{bt3}
\end{align}
From \eqref{four}, with the parameters $[1, 5, 9, 15, 33]$ and
with $q$ replaced by $q^{72}$, we find that
\begin{equation}
[4:72]\big\{[14,32,34:72]+q^4[6,24,30:72]\big\}=[8,10,28,34:72].\label{bt4}
\end{equation}
Returning \eqref{bt3} and substituting \eqref{bt4}, we arrive at
\begin{align}
&[8,10,11,17,19,25,28,32:72]\notag\\
&=[5,16,31:72]\big\{[8,13,23,28,34:72]+q^5[6,12,17,19,24:72]\big\}.\label{bt5}
\end{align}
By switching to base $q^{36}$, we have from \eqref{bt5} that
\begin{align*}
&[5,8,11,16,17:36](5,16:36)\notag\\
&=[5,8:36](8:36)\big\{[8,13,17:36](17:36)+q^5[3,12,17:36](3:36)\big\}.
\end{align*}
After cancellation, we are led to prove
\begin{equation}
[11,16:36](5,16:36)=(8:36)\big\{[8,13:36](17:36)+q^5[3,12:36](3:36)\big\}.\label{bitt}
\end{equation}
Now, \eqref{bitt} follows from \eqref{four} with $q$ replaced by
$q^{36}$ and $a$, $b$, $c$, $x$, $y$ replaced by $q^{17}$,
$-q^{20}$, $q^{4}$, $q$ $q^{28}$ $-q^{12}$, respectively.
\end{proof}

\section{Modulus $M = 80$}
There is only one class  of identities under multiplication by the
group $U(80)$.
\begin{theorem}\label{80.1}
\begin{align*}
\mathrm{(i)}\hspace{0.5in}&\text{Let}\; S\equiv\pm\{4, 5, 6, 7, 15, 17, 23, 25, 26, 33, 35, 36\}\pmod {80},\; \text{and} \\
&T\equiv\pm\{2, 5, 11, 12, 15, 18, 19, 21, 25, 28, 29, 35\} \pmod{80}.\; \text{Then}\\
&p(S,n)=p(T,n),\; \text{for all}\; n\neq 2.
\end{align*}
\begin{align*}
\mathrm{(ii)}\hspace{0.5in}&\text{Let}\; S\equiv\pm\{3, 4, 5, 13, 14, 15, 25, 27, 34, 35, 36, 37\}\pmod {80},\; \text{and} \\
&T\equiv\pm\{1, 5, 9, 12, 15, 22, 25, 28, 31, 35, 38, 39\} \pmod{80}.\; \text{Then}\\
&p(S,n)=p(T,n-3),\; \text{for all}\; n\geq 3.
\end{align*}
\end{theorem}
\begin{proof}
The identities (i) and (ii) follow from \eqref{four2} with the
choice of variables $[1, 5, 7, 12, 19]$, $[1, 5, 2, 7, 14]$ and with
$q$ replaced by $q^{40}$ in both instances.
\end{proof}

\section{Modulus $M = 82$}
There is only one class  of identities under multiplication by the
group $U(82)$.
\begin{theorem}\label{82.1}
\begin{align*}
\mathrm{(i)}\hspace{0.5in}&\text{Let}\; S\equiv\pm\{6, 8, 10, 11, 13, 15, 17, 28, 29, 31, 33, 35\}\pmod {82},\; \text{and} \\
&T\equiv\pm\{2, 5, 11, 17, 18, 23, 24, 30, 31, 33, 37, 39\} \pmod{82}.\; \text{Then}\\
&p(S,n)=p(T,n-6),\; \text{for all}\; n\geq 6.
\end{align*}
\begin{align*}
\mathrm{(ii)}\hspace{0.5in}&\text{Let}\; S\equiv\pm\{1, 3, 8, 9, 10, 14, 21, 25, 27, 31, 33, 38\}\pmod {82},\; \text{and} \\
&T\equiv\pm\{1, 3, 7, 9, 11, 17, 19, 24, 27, 30, 32, 40\} \pmod{82}.\; \text{Then}\\
&p(S,n)=p(T,n),\; \text{for all}\; n\neq 7.
\end{align*}
\begin{align*}
\mathrm{(iii)}\hspace{0.5in}&\text{Let}\; S\equiv\pm\{3, 4, 5, 13, 14, 15, 27, 29, 35, 36, 37, 38\}\pmod {82},\; \text{and} \\
&T\equiv\pm\{1, 5, 9, 12, 15, 23, 26, 29, 32, 37, 39, 40\} \pmod{82}.\; \text{Then}\\
&p(S,n)=p(T,n-3),\; \text{for all}\; n\geq 3.
\end{align*}
\begin{align*}
\mathrm{(iv)}\hspace{0.5in}&\text{Let}\; S\equiv\pm\{2, 3, 7, 13, 18, 19, 22, 23, 27, 34, 35, 39\}\pmod {82},\; \text{and} \\
&T\equiv\pm\{1, 6, 9, 13, 16, 20, 21, 23, 25, 28, 35, 39\}\pmod{82}.\; \text{Then}\\
&p(S,n)=p(T,n),\; \text{for all}\; n\neq 1.
\end{align*}
\begin{align*}
\mathrm{(v)}\hspace{0.5in}&\text{Let}\; S\equiv\pm\{7, 11, 12, 15, 16, 17, 19, 20, 21, 25, 26, 29\}\pmod {82},\; \text{and} \\
&T\equiv\pm\{4, 5, 7, 19, 21, 22, 25, 31, 33, 34, 36, 37\}\pmod{82}.\; \text{Then}\\
&p(S,n)=p(T,n-7),\; \text{for all}\; n\geq 7.
\end{align*}

\end{theorem}
\begin{proof}
The identities (i)--(v) follow from \eqref{four2} with the choice
of variables $[1, 9, 3, 13, 14]$, $[1, 2, 9, 12, 18]$, $[1, 5, 2, 7,
14]$, $[1, 3, 4, 10, 17]$, $[1, 8, 12, 32, 34]$,  and with $q$
replaced by $q^{41}$ in each instance.
\end{proof}

{\textbf{Acknowledgement:}}
We would like to thank the referee for a careful reading of the
paper and for a number of suggestions.

\end{document}